\documentclass[10pt]{amsart}
\usepackage{calrsfs}

\usepackage{amsmath,tikz}
\usetikzlibrary{tikzmark,fit}
\newcommand\bigzero{\makebox(0,0){\text{\huge0}}}

\usepackage{graphicx}
\usepackage{amsthm}
\usepackage{amssymb}
\usepackage{multirow}
\usepackage{tikz-cd}
\usepackage{amsmath}
\usepackage{todonotes}
\usepackage{amsbsy}
\usepackage[all]{xy}
\usepackage{qtree}

\usepackage{enumitem}
\usepackage{xfrac}    
\usepackage{color, colortbl}
\definecolor{LightCyan}{rgb}{0.88,1,1}
\definecolor{Gray}{gray}{0.9}

\usepackage{faktor}

\usepackage{changes}
\definechangesauthor[color=orange,name={Aristidesd Kontogeorgis}]{AK}
\definechangesauthor[color=blue,name={Alex Terezakis}]{AT}

\newtheorem{theorem}{Theorem}
\newtheorem{lemma}[theorem]{Lemma}

\newtheorem{proposition}[theorem]{Proposition}
\theoremstyle{definition}

\newtheorem{remark}[theorem]{Remark}
\newtheorem{definition}[theorem]{Definition}

\newcommand{\Z}{\mathbb{Z}}

\newcommand{\asp}{ \begin{array}{c} \; \\ \; \\mathrm{end}{array}} 

\renewcommand{\mod}{{\;\mathrm{mod}}}

\newcommand\bovermat[2]{%
  \makebox[0pt][l]{$\smash{\overbrace{\phantom{%
    \begin{matrix}#2\end{matrix}}}^{\text{#1}}}$}#2}

\date{\today}

\title{On the lifting problem of representations of a metacyclic group.}

\author[A. Kontogeorgis]{Aristides Kontogeorgis}
\address{Department of Mathematics, National and Kapodistrian  University of Athens
Pane\-pist\-imioupolis, 15784 Athens, Greece}
\email{kontogar@math.uoa.gr}

\author[A. Terezakis]{Alexios Terezakis }
\address{Department of Mathematics, National and Kapodistrian University of Athens\\
Panepistimioupolis, 15784 Athens, Greece}
\email{aleksistere@math.uoa.gr}

\date \today

\makeatletter
\newcommand{\aprod}{\mathop{\operator@font \hbox{\Large$\ast$}}}
\makeatother

\begin{document}


\begin{abstract}
We give a necessary and sufficient condition for a modular representation of a group $G=C_{p^h} \rtimes C_m$ in a field of characteristic zero to be lifted to a representation over local principal ideal domain of characteristic zero containing the $p^h$ roots of unity. 
\end{abstract}
\maketitle

\section{Introduction}
The lifting problem for a representation 
\[
  \rho: G \rightarrow \mathrm{GL}_n(k), 
\]
where $k$ is a field of characteristic $p>0$,
is about finding a local ring $R$ of characteristic $0$, with maximal ideal $\mathfrak{m}_R$ such that $R/\mathfrak{m}_R=k$, so that the following diagram is commutative:
\[
  \xymatrix{
   & \mathrm{GL}_n(R) \ar[d] \\
   G \ar[r] \ar[ru] & \mathrm{GL}_n(k)
  }
\]
Equivalently one asks if there is a free $R$-module $V$, which is also an $R[G]$-module such that $V\otimes_R R/\mathfrak{m}_R$ is the $k[G]$-module corresponding to our initial representation. 
We know that projective $k[G]$-modules lift in characteristic zero, \cite[chap. 15]{SerreLinear}, but for a general $k[G]$-module 
such a lifting  is not always possible, for example, see \cite[prop. 15]{https://doi.org/10.48550/arxiv.2101.11084}. This article aims to study the lifting problem for the group $G=C_q \rtimes C_m$, where $C_q$ is a cyclic group of order $p^h$ and $C_m$ is a cyclic group of order $m$, $(p,m)=1$ and give necessary and sufficient condition in order to lift. We assume that the local ring $R$ contains the $q$-roots of unity and $k$ is algebraically closed, and we might need to consider a ramified extension of $R$, in order to ensure that certain $q$-roots of unit are distant in the $\mathfrak{m}_R$-topology, see remark 
\ref{cond-ram-R}. An example of such a ring $R$ is the ring of Witt vectors $W(k)[\zeta_q]$ with the $q$-roots of unity adjoined to it. 

We notice that a decomposable $R[G]$-module $V$ gives rise to a decomposable $R$-module modulo 
$\mathfrak{m}_R$ and also an indecomposable
 $R[G]$-module can break in the reduction modulo $\mathfrak{m}_R$
into a direct sum of indecomposable $k[G]$-summands.  We also give a classification of $k[C_q \rtimes C_m]$-modules in terms of Jordan decomposition and give the relation with the more usual uniserial description in terms of their socle \cite{alperin}. 

Our interest to this problem comes from the problem of lifting local actions. The local lifting problem considers the following question:
Does there exist an extension $\Lambda/W(k)$, and a representation 
\[
\tilde{\rho}: G \hookrightarrow \mathrm{Aut}(\Lambda[[T]]),
\]
such that if $t$ is the reduction of $T$, then the action of $G$ on $\Lambda[[T]]$
reduces to the action of $G$ on $k[[t]]$?

If the answer to the above question is positive, then we say that the $G$-action lifts to characteristic zero. 
A group $G$ for which every local $G$-action on $k[[t]]$ lifts to characteristic zero is called {\em a local Oort group for $k$}. Notice that cyclic groups are always local Oort groups. This result was known as the ``Oort conjecture'', which was recently proved by F. Pop \cite{MR3194816} using the work of A. Obus and S. Wewers \cite{ObusWewers}.

There are a lot of obstructions that prevent a local action to lift in characteristic zero. 
Probably the most important of these obstructions in the KGB-obstruction \cite{MR2919977}. It is believed that this is the only obstruction for the local lifting problem, see \cite{Obus12},
 \cite{Obus2017}. 
 {\em In \cite[Thm. 3]{https://doi.org/10.48550/arxiv.2101.11084} the authors have given a criterion for the local lifting which involves the lifting of a linear representation of the same group. }
 The case $G=C_q \rtimes C_m$ and especially the case of dihedral groups $D_q=C_q \rtimes C_2$, is a problem of current interest in the theory of local liftings, see \cite{Obus2017}, \cite{1912.12797}, \cite{MR3874854}.
For more details on the local lifting problem we refer to \cite{MR2441248}, \cite{MR2919977}, \cite{MR3591155}, \cite{Obus12}.  

Keep also in mind that the $C_q \rtimes C_m$  groups were important to the study of group actions in holomorphic differentials of curves defined over fields of positive characteristic $p$, where the group involved has cyclic $p$-Sylow subgroup, see \cite{MR4130074}.

Let us now describe the method of proof. 
For understanding the splitting of indecomposable 
$R[G]$-modules modulo $\mathfrak{m}_R$, we develop a version of Jordan normal form in lemma \ref{lemma:Jordan}  for endomorphisms $T:V \rightarrow V$ of order $p^h$, where $V$ is a free module of rank $d$. We give a way to select this basis, by selecting an initial suitable  element $E \in V$, see lemma \ref{lemma:lin-inde}. The normal form (as given in eq. (\ref{eq-matrix-lift})) of the element $T$ of order $q$,  determines the decomposition of the reduction. We show that for every indecomposable summand $V_i$ of $V$, we can select $E$ as an eigenvalue of the generator $\sigma$ of $C_m$ and then by forcing the  relation $\Gamma T =T^\alpha \Gamma$ to hold, we see how the action of $\sigma$ can be extended recursivelly to an action of $\sigma$ on $V_i$, this is done in lemma \ref{lemma:recursive}. 
Proving that this gives indeed a well defined action is a technical computation and is done in lemmata \ref{lemma:sum-prod}, \ref{lemma:ABlL}, \ref{zeroGamma}, \ref{d-1}, \ref{lemma:30}.

The important thing here is that the definition of the action of $\sigma$ on $E$ is the ``initial condition'' of a dynamical system that determines the action of $C_m$ on the indecomposable summand $V_i$. The $R[C_q \rtimes C_m]$ indecomposable module $V_i$ can break into a direct sum $V_{\alpha}(\epsilon_\nu,\kappa_\nu)$-modules $1\leq \nu \leq s$ (for a precise definition of them see definition \ref{def:V}, notice that $\kappa_i$ denotes the dimension). The action of $\sigma$ on each $V_{\alpha}(\epsilon_\nu,\kappa_\nu)$ can be uniquely determined by the action of $\sigma$ on an initial basis element as shown in section \ref{Vdesc}, again by a ``dynamical system'' approach,  where we need $s$ initial conditions, one for each $V_\alpha(\epsilon_\nu,\kappa_\nu)$.  The lifting condition essentially means that the indecomposable summands 
$V_\alpha(\epsilon,\kappa)$ of the special fibre, should be able to be rearranged in a suitable way, so that they can be obtained as reductions of indecomposable $R[C_q \rtimes C_q]$-modules.
The precise expression of our lifting criterion is given in the following proposition: 
 \begin{proposition}
\label{cond-lift}
Consider a $k[G]$-module $M$ which is decomposed as a direct sum 
\[
  M= V_{\alpha}(\epsilon_1,\kappa_1) \oplus \cdots \oplus V_{\alpha}(\epsilon_s,\kappa_s). 
\]
The module lifts to an $R[G]$-module if and only if 
the set $\{1,\ldots,s\}$ can be written as a disjoint union of sets $I_\nu$, $1\leq \nu \leq t$ so that

\begin{enumerate}[label=\alph*]
\item  \label{cond-lifta}$\sum_{\mu \in I_\nu} \kappa_\mu \leq q$, for all $1\leq \nu \leq t$. 
\item \label{cond-liftb}$\sum_{\mu \in I_\nu} \kappa_\mu\equiv a \mod m$  for all $1\leq \nu \leq t$, where $a\in \{0,1\}$. 
\item \label{cond-liftc} 
For each $\nu$, $1\leq \nu \leq t$ there is an enumeration $\sigma: \{1,\ldots,\#I_\nu\} \rightarrow I_\nu\subset \{1,..,s\}$,
such that 
\[
  \epsilon_{\sigma(2)}= \epsilon_{\sigma(1)} \alpha^{\kappa_{\sigma(1)}},
  \epsilon_{\sigma(3)}= \epsilon_{\sigma(3)} \alpha^{\kappa_{\sigma(3)}},\ldots
  \epsilon_{\sigma(s)}= \epsilon_{\sigma(s-1)} \alpha^{\kappa_{\sigma(s-1)}}
\]
\end{enumerate}
\end{proposition}
In the above proposition, each set $I_\nu$ corresponds to a collection of modules $V_{\alpha}(\epsilon_\mu,\kappa_\mu)$, $\mu \in I_\nu$ which come as the reduction of an indecomposable  $R[C_q \rtimes C_m]$-module $V_\nu$ of $V$. 

\noindent {\bf Aknowledgements} A. Terezakis is a recipient of financial support in the context of a doctoral thesis (grant number MIS-5113934). The implementation of the doctoral thesis was co-financed by Greece and the European Union (European Social Fund-ESF) through the Operational Programme—Human Resources Development, Education and Lifelong Learning—in the context of the Act—Enhancing Human Resources Research Potential by undertaking a Doctoral Research—Sub-action 2: IKY Scholarship Programme for Ph.D. candidates in the Greek Universities.

\bigskip
\begin{center}
\includegraphics[scale=0.23]{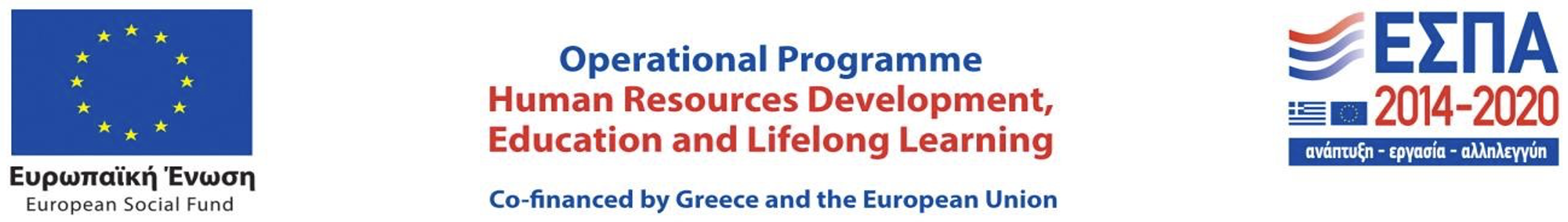}
\end{center}
\section{Notation}
 Let $\tau$ be a generator of the cyclic group $C_q$ and $\sigma$ be a generator of the cyclic group $C_m$.  The group $G$ is given in terms of generators and relations as follows:
\[
G=\langle \sigma,\tau | \tau^q=1, \sigma^m=1, \sigma \tau \sigma^{-1}= \tau^{\alpha} \mbox{ for some } \alpha\in \mathbb{N},
 1\leq \alpha \leq p^h-1, (\alpha, p)=1 \rangle.
\]
The integer $\alpha$ satisfies the following congruence:
\begin{equation}
  \label{a-modq}
  \alpha^m \equiv 1 \mod q
\end{equation}
as one sees by computing $\tau=\sigma^m \tau\sigma^{-m}=\tau^{\alpha^m}$. 
Also the integer $\alpha$ can be seen as an element in the finite field $\mathbb{F}_p$, and it is a
 $(p-1)$-th root of unity, not necessarily primitive. In particular the following holds:
\begin{lemma}  \label{loga} Let $\zeta_m\in k$ be a fixed primitive $m$-th root of unity. 
There is a natural number  $a_0$, $0\leq a_0 < m-1$ such that $\alpha=\zeta_m^{a_0}$.  
\end{lemma}
\begin{proof}
The integer $\alpha$ if we see it as an element in $k$ is an element in the finite field $\mathbb{F}_p\subset k$, therefore $\alpha^{p-1}=1$ as an element in $\mathbb{F}_p$.
Let $\mathrm{ord}_p(\alpha)$ be the order of $\alpha$ in $\mathbb{F}_p^*$. By eq. (\ref{a-modq}) we have that $\mathrm{ord}_p(\alpha) \mid p-1$ and $\mathrm{ord}_p(\alpha) \mid m$, that is $\mathrm{ord}_p(\alpha) \mid (p-1,m)$. 

The primitive $m$-th root of unity $\zeta_m$ generates a finite field $\mathbb{F}_p(\zeta_m)=\mathbb{F}_{p^\nu}$
for some integer $\nu$, which has cyclic multiplicative  group $\mathbb{F}_{p^\nu}\backslash \{0\}$ containing both the cyclic groups $\langle \zeta_m \rangle$ and $\langle \alpha \rangle$. Since for every divisor $\delta$ of the order of a cyclic group $C$ there is a unique subgroup $C' < C$ of order $\delta$ we have that $\alpha \in \langle \zeta_m \rangle$, and the result follows. 
\end{proof}

\begin{definition}
For each $p^i \mid q$ we define
$
\mathrm{ord}_{p^i} \alpha
$
to be the smallest natural number $o$ such that $\alpha^o\equiv 1 \mod p^i$.
\end{definition}
 It is clear that for $\nu \in \mathbb{N }$
\[
\alpha^{\nu} \equiv 1 \mod p^i \Rightarrow \alpha^\nu \equiv 1 \mod p^j \text{ for all } j \leq i. 
\]
Therefore
\[
\mathrm{ord}_{p^j}\alpha \mid \mathrm{ord}_{p^i} \alpha \text{ for } j \leq i.
\]
On the other hand $\alpha \in \mathbb{N}$ and $\alpha^{p-1} \equiv 1 \mod p$ so $\mathrm{ord}_p \alpha \mid p-1$. 
Also since $\sigma^t \tau \sigma^{-t}= \tau^{\alpha^t}$ we have that $\alpha^m \equiv 1 \mod p^h$, therefore
$\mathrm{ord}_p \alpha \mid \mathrm{ord}_{p^i}\alpha \mid \mathrm{ord}_{p^h} \alpha \mid m$, for $1\leq i \leq h$. 

\begin{lemma}
\label{lemma:center}
The center $\mathrm{Cent}_G(\tau)=\langle \tau, \sigma^{\mathrm{ord}_{p^h}\alpha}\rangle$. Moreover 
\[
\frac{
|\mathrm{Cent}_G(\tau)|}{p^h}
= \frac{m}{
\mathrm{ord}_{p^h}( \alpha)
}=:m'
\]
\end{lemma}
\begin{proof}
The result follows by observing
 $(\tau^\nu\sigma^t) \tau (\tau^\nu\sigma^{t})^{-1}=\tau^{\alpha^t}$, for all $1\leq \nu \leq q$, $1\leq t \leq m$. 
\end{proof}

\begin{remark}
If $\mathrm{ord}_p \alpha=m$ then  $\mathrm{ord}_{p^i} \alpha=m$ for all $1\leq i \leq h$.
\end{remark}

\begin{lemma}
If the group $G=C_q\rtimes C_m$ is a subgroup of $\mathrm{Aut}(k[ [t]])$, then all orders $\mathrm{ord}_{p^i} \alpha=m/m'$, for all $1\leq i \leq h$. 
\end{lemma}
\begin{proof}
We will use the notation of the book of J.P.Serre on local fields \cite{SeL}. By \cite[Th.1.1b]{MR2577662}
we have that the first gap $i_0$ in the lower ramification filtration of the cyclic group $C_q$ satisfies $(m,i_0)=m'$. 

The ramification relation \cite[prop. 9 p. 69]{SeL}
\[
\alpha \theta_{i_0}(\tau)=\theta_{i_0}(\tau^\alpha)=\theta_{i_0}(\sigma \tau \sigma^{-1})=\theta_0(\sigma)^{i_0}\theta_{i_0}(\tau),
\]
implies that $\theta_{0}(\sigma)^{i_0}=\alpha\in \mathbb{N}$.  
From $(m,i_0)=m'$ and the fact that $\mathrm{ord} \theta_0(\sigma)=m$ we obtain
\[
\frac{m}{m'}=\mathrm{ord} \theta_{0}(\sigma)^{i_0}=\mathrm{ord}_p(\alpha). 
\]
Thus 
\[
\frac{m}{m'} =\mathrm{ord}_{p}\alpha |  \mathrm{ord}_{p^i}\alpha | \mathrm{ord}_{p^h}\alpha =\frac{m}{m'}.  
\]
Hence all orders
$\mathrm{ord}_{p^i}\alpha=m/m'$. 
\end{proof}
\begin{remark}
If the KGB-obstruction vanishes and $\alpha\neq 1$, then by \cite{Obus12}[prop. 5.9] $i_0 \equiv -1 \mod m$ and  $\mathrm{ord}_{p^i}\alpha=m$ for all $1\leq i \leq h$. 
\end{remark}

\section{Indecomposable $C_q\rtimes C_m$ modules, modular representation theory}

\label{Vdesc}
In this section we will describe the indecomposable $C_q\rtimes C_m$-modules. We will give two methods in studying them. The first one is needed since it is in accordance to the method we will give in order to describe indecomposable $R[C_q\rtimes C_m]$-modules. The second one, using the structure of the socle, is the standard method of describing $k[C_q \rtimes C_m]$-modules in modular representation theory.  

\subsection{Linear algebra method.}

The indecomposable modules of the $C_q$ are determined by the Jordan normal forms of the generator $\tau$
of the cyclic group $C_q$. So for each $1\leq \kappa \leq p^h$ there is exactly one $C_q$ indecomposable module denoted by $J_\kappa$. Therefore we have the following decomposition of an indecomposable $C_q \rtimes C_m$-module $M$ considered as a $C_q$-module. 
\begin{equation} \label{decMa}
M= J_{\kappa_1} \oplus \cdots \oplus J_{\kappa_r}. 
\end{equation}
\begin{lemma}
In the indecomposable module $J_\kappa$ for every element $E$ such that
 \[
 (\tau- \mathrm{Id}_{\kappa_i})^{\kappa_i-1}E \neq 0
 \]
the elements $B=\{E, (\tau-\mathrm{Id}_\kappa)E,\ldots, (\tau- \mathrm{Id}_{\kappa})^{\kappa-1}E\}$ form a basis of $J_\kappa$ such that the matrix of $\tau$ with respect to this basis is given by 
\begin{equation}
\label{eq:tau-form}
\tau=\mathrm{Id}_\kappa+
\begin{pmatrix}
0 & \cdots &   \cdots & \cdots &  0 \\
1 & \ddots  &   & & \vdots  \\ 
0 &  \ddots  & \ddots   &  & \vdots \\
\vdots & \ddots &   1 & 0 & \vdots \\
0 & \cdots & 0 & 1 & 0
\end{pmatrix}.
\end{equation}
\end{lemma}
\begin{proof}
Since the set $B$ has $k$-elements it is enough to prove that it consists of linear independent elements. Indeed, consider a linear relation
\[
  \lambda_0 E +\lambda_1(\tau-\mathrm{Id}_\kappa)E + \cdots + \lambda_{\kappa-1}(\tau-\mathrm{Id}_\kappa)^{\kappa-1}E =0.\]
  By applying $(\tau-\mathrm{Id}_\kappa)^{\kappa-1}$ we obtain $\lambda_0 (\tau-\mathrm{Id}_\kappa)^{\kappa-1}=0$, which gives us $\lambda_0=0$. We then apply $(\tau-\mathrm{Id}_\kappa)^{\kappa-2}$ to the linear relation and by the same argument we obtain $\lambda_1=0$ and we continue this way proving that $\lambda_0=\cdots=\lambda_{\kappa-1}=0$. The matrix form of $\tau$ in this basis is immediate.  
\end{proof}
We will now prove that $\sigma$ acts on each $J_\kappa$ of eq. (\ref{decMa}) proving that $r=1$. Since the field $k$ is algebraically closed and $(m,p)=1$ we know that there is a basis of $M$ consisting of eigenvectors of $\sigma$.
There is an eigenvector $E$ of $\sigma$, which is not in the kernel of $(\tau-\mathrm{Id}_\kappa)^{\kappa_1-1}$. Then 
the elements of the set $B=\{E,(\tau-\mathrm{Id}_\kappa)E,\ldots,(\tau-\mathrm{Id}_\kappa)^{\kappa_1-1}E\}$ are linearly independent and form a direct $C_q$ summand of $M$ isomorphic to $J_{\kappa_1}$. 

We will now show that this module is an $k[C_q \rtimes C_m]$-module. For this, we have to show that the generator $\sigma$ of $C_m$ acts on the basis $B$.  
Observe that for every $0\leq i \leq \kappa_1-1 <p^h$
\[
\sigma (\tau-1)^{i-1}= (\tau^\alpha-1)^{i-1} \sigma.
\]
Set $e=E_1$ and $\kappa=\kappa_1$.
This means that the action of $\sigma$ on $e$ determines the 
action of $\sigma$ on all other basis elements $e_\nu:=(\tau-1)^{\nu-1}e$, 
$1\leq \nu \leq \kappa_1$.

Let us compute:
\[
\sigma e_{i+1} =  \sigma(\tau-1)^ie=(\tau^\alpha-1)^i \zeta_m^\lambda e
\]
On the basis $\{e_1,\ldots,e_{\kappa_1}\}$ the matrix $\tau$ is given by eq. (\ref{eq:tau-form})
hence using the binomial formula we compute 
\begin{equation}
\label{eq:a-power-t-spec}
\tau^{\alpha}=
\begin{pmatrix}
1 & 0 &  \cdots &  \cdots  & \cdots &  0 \\
 \binom{\alpha}{1}& 1 & \ddots  &  &  & \vdots \\
\binom{\alpha}{2} & \binom{\alpha}{1} & \ddots & \ddots &  & \vdots \\
\binom{\alpha}{3} & \binom{\alpha}{2} & \ddots & 1 &   \ddots & \vdots\\
\vdots & \vdots  & \ddots & \binom{\alpha}{1} & 1  & 0 \\
 \binom{\alpha}{k}& \binom{\alpha}{k-1} & \cdots & \binom{\alpha}{2} & \binom{\alpha}{1} & 1
\end{pmatrix}.
\end{equation}
Thus $\tau^\alpha-1$ is a nilpotent matrix $A=(a_{ij})$ of the form:
\[
a_{ij}=
\begin{cases}
\binom{\alpha}{\mu} & \mbox{ if } j= i-\mu \text{ for some } \mu, 1\leq \mu \leq \kappa \\
0 & \mbox{ if } j \geq i
\end{cases}
\]
The $\ell$-th power $A^\ell=(a_{ij}^{(\ell)})$ of $A$ is then computed by (keep in mind that $a_{ij}=0$ for $i\leq j$)
\[
a_{ij}^{(\ell)}=\sum_{i < \nu_1 < \cdots < \nu_{\ell-1} < j} 
a_{i,\nu_1} a_{\nu_1,\nu_2} a_{\nu_2,\nu_3} \cdots a_{\nu_{\ell-1}, j} 
\] 
This means that $i-j>\ell$ in order to have $a_{ij}\neq 0$. 
Moreover for $i=j+\ell$ (which is the the first non zero diagonal below  the main diagonal) we have 
\[
a_{i,i+\ell}= a_{i,i+1} a_{i+1,i+2} \cdots a_{i+\ell-1, i+\ell}=\binom{\alpha}{1}^\ell=\alpha^\ell.
\]
Therefore, the matrix of $A^{\ell}$ is of the following form:
\vskip 0.2cm
\begin{equation} \label{powerA}
\begin{pmatrix}
\bovermat{$k-\ell$}{0        &\cdots & \cdots & 0} & \bovermat{$\ell$}{0 & \cdots & 0} \\
\vdots   &       &       & \vdots & \vdots & & \vdots\\ 
0        &\cdots & \cdots & 0 & 0&  \cdots & 0 \\
\alpha^\ell & \ddots & & 0 &\vdots & &  \vdots \\
* & \alpha^\ell & \ddots & \vdots &\vdots & & \vdots\\
\vdots & \ddots & \ddots & 0 & \vdots &   & \vdots \\
* & \cdots & * & \alpha^\ell & 0 & \cdots & 0\\
\end{pmatrix}
\end{equation}
\begin{definition}
\label{def:V}
We will denote by $V_\alpha(\lambda,\kappa)$ the indecomposable $\kappa$-dimensional $G$-module given by the basis elements
$\{(\tau-1)^\nu e, \nu=0,\ldots,\kappa-1\}$, where $\sigma e=\zeta_m^\lambda  e$. 
\end{definition}

This definition is close to 
the notation used in \cite{SK}.

\begin{lemma} \label{action-tau}
The action of $\sigma$ on the basis element $e_i$ of $V_\alpha(\lambda,\kappa)$ is given by:
\begin{equation} \label{actionCom}
\sigma e_i= \alpha^{i-1} \zeta_m^\lambda e_i + \sum_{\nu=i+1}^\kappa a_\nu e_\nu,   
\end{equation}
for some coefficients $a_i \in k$.
In particular the matrix of $\sigma$ with respect to the basis $e_1,\ldots,e_\kappa$ is lower triangular.
\end{lemma}
\begin{proof}
Recall that $e_i=(\tau-1)^{i-1} e_1$. Therefore 
\[
\sigma e_i=\sigma (\tau-1)^{i-1} e_1=(\tau^\alpha-1)^{i-1} \sigma e_1= \zeta_m^\lambda (\tau^\alpha-1)^{i-1} e_1. 
\]
The result follows by eq. (\ref{powerA})
\end{proof}

We have constructed a set of indecomposable modules $V_\alpha(\lambda,\kappa)$. Apparently $V_\alpha(\lambda,\kappa)$ can not be isomorphic to $V_\alpha(\lambda',\kappa')$ if $\kappa\neq \kappa'$, since they have different dimensions. 

Assume now that $\kappa=\kappa'$. Can the modules $V_\alpha(\lambda,\kappa)$ and $V_\alpha (\lambda',\kappa)$ be isomorphic for $\lambda \neq \lambda'$?

The eigenvalues of the prime to $p$ generator $\sigma$ on $V_\alpha(\lambda,\kappa)$are 
\[\zeta_m^{\lambda},\alpha\zeta_m^{\lambda},
\ldots, \alpha^{\kappa-1}\zeta_m^\lambda.\]
Similarly the  eigenvalues for $\sigma$ when acting on $V_\alpha(\lambda',\kappa)$ are 
\[
\zeta_m^{\lambda'},\alpha\zeta_m^{\lambda'},
\ldots, \alpha^{\kappa-1}\zeta_m^{\lambda'}.
\]
If the two sets of eigenvalues are different then the modules can not be isomorphic. But even if $\lambda \neq \lambda' \mod n$ the 
two sets of eigenvalues can still be equal. Even in this case the modules can 
not be isomorphic.

\begin{lemma}
The modules $V_\alpha(\lambda_1,\kappa)$ and $V_\alpha(\lambda_2,\kappa)$
are isomorphic if and only if $\lambda_1\equiv \lambda_2 \mod m$. 
\end{lemma}
\begin{proof}
Indeed, the module $V_\alpha(\lambda_1,\kappa)$ has an eigenvector for the 
action of $\sigma$ which generates the $V_\alpha(\lambda_1,\kappa)$ by powers of $(\tau-1)$, i.e. the vectors
\begin{equation} \label{SigmaBasis}
e, (\tau-1)e, (\tau-1)^2 e, \ldots, (\tau-1)^{\kappa-1}e 
\end{equation}
form a basis of $V_\alpha(\lambda_1,\kappa)$.

The elements $E$ which can generate $V_\alpha(\lambda_1,\kappa)$
by powers of $(\tau-1)$ are linear combinations 
\[
E =\sum_{\nu=0}^{\kappa-1} \lambda_i (\tau-1)^\nu e, 
\]
for $\lambda_i \in k$ and $\lambda_0\neq 0$. 

On the other hand  using eq. (\ref{actionCom}) we see that $\sigma$ with respect to the basis given in eq. (\ref{SigmaBasis})
admits the matrix form:
\[
\begin{pmatrix}
\zeta_m^{\lambda} & 0 & \cdots & \cdots &  0  \\
0  & \alpha \zeta_m^{\lambda}  & 0 & \cdots &  0 \\ 
\vdots & \ddots & \ddots & \ddots & \vdots \\
\vdots &  & \ddots  &  \ddots & \vdots \\
0  & \cdots & \cdots & 0 & \alpha^{\kappa-1} \zeta_m^{\lambda}
\end{pmatrix}.
\]
It is now easy to see from the above matrix that 
every eigenvector of the 
eigenvalue $\alpha^\nu \lambda_1$, $\nu>1$ is expressed as a linear 
combination of the basis given in eq. (\ref{SigmaBasis}), where the coefficient of $e$ is  zero. 

Therefore, the eigenvector of the eigenvalue $\alpha^\nu \zeta_m$ 
can not generate the module $V_\alpha(\lambda,\kappa)$ by powers of 
$(\sigma-1)^\nu$. 
\end{proof}

\subsection{The uniserial description}
We will now give an alternative description of the indecomposable $C_q \rtimes C_m$-modules, which is used in  
\cite{MR4130074}. 

It is known that $\mathrm{Aut}(C_q)\cong \mathbb{F}_p^* \times Q$, for some abelian $p$-group $Q$. The representation $\psi:C_m \rightarrow \mathrm{Aut}(C_q)$ given by the action of $C_m$ on $C_q$ is known to factor through a character $\chi:C_m \rightarrow \mathbb{F}_p^*$. The order of $\chi$
divides $p-1$ and $\chi^{p-1}=\chi^{-(p-1)}$ is the trivial one dimensional character.

For all $i\in \Z$, $\chi^i$ defines a simple $k[C_m]$-module of $k$ dimension one, which we will denote by $S_{\chi^i}$. For $0 \leq \ell \leq m-1$ denote by $S_\ell$ the simple module where on which $\sigma$ acts as $\zeta_m^\ell$.  Both $S_{\chi^i}$, $S_\ell$ can be seen as $k[C_q\rtimes C_m]$-modules using inflation. Finally for $0\leq \ell \leq m-1$ we define $\chi^i(\ell) \in \{0,1,\ldots,m-1\}$ such that $S_{\chi^i(\ell)} \cong S_\ell \otimes_k S_{\chi^i}$. 

There are $q\cdot m$ isomorphism classes of indecomposable $k[C_q \rtimes C_m]$-modules and are all uniserial. An indecomposable $k[C_q \rtimes C_m]$-module $U$ is unique determined by its socle, which is the kernel of the action of $\tau-1$ on $U$, and its $k$-dimension. For $0 \leq \ell \leq m-1$ and $1\leq \mu \leq q$, let $U_{\ell,\mu}$ be the indecomposable $k[C_q\rtimes C_m]$
module with socle $S_a$ and $k$-dimension $\mu$. Then $U_{\ell,\mu}$ is uniserial and its $\mu$
ascending composition factors are the first $\mu$ composition factors of the sequence 
\[
S_\ell, S_{\chi^{-1}(\ell)},S_{\chi^{-2}(\ell)},\ldots,S_{\chi^{-(p-2)}(\ell)}, 
S_{\ell}, S_{\chi^{-1}(\ell)},S_{\chi^{-2}(\ell)},\ldots, S_{\chi^{-(p-2)}(\ell)}.
\]
Notice that in our notation $V_{\alpha}(\lambda,\kappa)=U_{\lambda+\kappa,\kappa}$.

\begin{remark}
The condition $\mathrm{ord}_{p^i}=m$ for all $1\leq i \leq h$, is equivalent to requiring that $\psi_i:C_m \rightarrow \mathrm{Aut}(C_{p^i})$ is faithful for all $i$.   
\end{remark}

\section{Lifting of representations}
\label{sec:liftrep}
\begin{proposition}
\label{prop:div-obstr}
Let $G=C_q \rtimes C_m$. Assume that for all $1\leq i \leq h$, $ord_{p^i}a=m$. 
If the $G$-module $V$ lifts to  an $R[G]$-module $\tilde{V}$, where $K=\mathrm{Quot(R)}$ is a field of characterstic zero, then 
\[
m \mid 
\left(
\dim (\tilde{V}\otimes_R K) - \dim (\tilde{V}\otimes_R K)^{C_q}
\right).
\]
Moreover, if $\tilde{V}(\zeta_q^{\alpha^i\kappa})$ is the eigenspace  of the 
eigenvalue $\zeta_q^{\alpha^i\kappa}$ of $T$ acting on $\tilde{V}$, then 
\[
\dim \tilde{V}(\zeta_q^\kappa)=
\dim \tilde{V}(\zeta_q^{\alpha\kappa})=
\dim \tilde{V}(\zeta_q^{\alpha^2\kappa})=
\cdots =
\dim \tilde{V}(\zeta_q^{\alpha^{m-1}\kappa}).
\]  
\end{proposition}
\begin{proof}
Consider a lifting $\tilde{V}$ of $V$. The generator $\tau$ of the cyclic part $C_q$ has  eigenvalues  $\lambda_1,\ldots,\lambda_s$ which are $p^n$-roots of unity. Let $\zeta_q$ be a primitive $q$-root of unity. 
Consider any eigenvalue $\lambda \neq 1$. It is of the form $\lambda=\zeta_q^{\kappa}$ for some $\kappa \in \mathbb{N}, q\nmid \kappa$. 
If $E$ is an eigenvector of $T$ corresponding to $\lambda$, that is $\tau E= \zeta_q^{\kappa} E$ then 
\[
\tau \sigma^{-1} E= \sigma^{-1} \tau^\alpha E= \zeta_q^{\kappa \alpha^{m-1}} \sigma^{-1}E 
\]
and we have a series of eigenvectors $E, \sigma^{-1} E, \sigma^{-2} E, \cdots$
with corresponding eigenvalues $\zeta_q^\kappa, \zeta_q^{\kappa \alpha}, \zeta_q^{\kappa a^2} \cdots, \zeta_q^{ \kappa \alpha^o}$, where 
$o=\mathrm{ord}_{q/(q,k)} $. 
Indeed, the integer $o$ satisfies the
\[
\kappa \alpha^o \equiv \kappa \mod q \Rightarrow \alpha^m \equiv 1 \mod \frac{q}{(q,k)}.
\]
Therefore the eigenvalues $\lambda\neq 1$ form orbits of size $m$, while the eigenspace of the eigenvalue $1$ is just the invariant space $V^G$ and the result follows.
\end{proof}

\section{Indecomposable $C_q\rtimes C_m$ modules, integral representation theory}
From now on $V$ be a free $R$-module, where $R$ is an integral local principal ideal domain with maximal ideal $\mathfrak{m}_R$, $R$ has characteristic zero and that $R$ contains all $q$-th roots of unity and has characteristic zero.  Let $K=\mathrm{Quot}(R)$. 

The indecomposable modules for a cyclic group both in the ordinary and in the modular case are described by writing down the Jordan normal form of a generator of the cyclic group. Since in integral representation theory there are infinitely many non-isomorphic indecomposable $C_q$-modules for $q=p^h$, $h\geq 3$,  one is not expecting to have a theory of Jordan normal forms even if one works over complete local principal ideal domains \cite{MR0140575}, \cite{MR0144980}. 

\begin{lemma}
Let $T$ be an element of order $q=p^h$ in $\mathrm{End}(V)$, then the minimal polynomial of $T$ has simple eigenvalues and $T$ is diagonalizable when seen as an element in $\mathrm{End}(V\otimes K)$.
\end{lemma}
\begin{proof}
Since $T^q=\mathrm{Id}_V$, the minimal polynomial of $T$ divides $x^q-1$, which has simple roots over a field of characteristic zero. This ensures that $T\in \mathrm{End}(V\otimes K)$ is diagonalizable. 
\end{proof}

\begin{lemma}
\label{lemma:lin-inde}
Let $f(x)=(x- \lambda_1)(x- \lambda_2)\cdots (x - \lambda_d)$ be the minimal polynomial of $T$ on $V$. 
There is an element $E\in V$, such that 
\[
E, (T- \lambda_1\mathrm{Id}_V)E, (T- \lambda_2\mathrm{Id}_V)(T- \lambda_1\mathrm{Id}_V)E, \ldots, (T- \lambda_{d-1}\mathrm{Id}_V)\cdots (T- \lambda_1 \mathrm{Id}_V)E
\]
 are linear independent elements in $V\otimes K$. 
\end{lemma}
\begin{proof}
Consider the endomorphisms for $i=1,\ldots,d$
\[
\Pi_i=\prod_{\nu=1 \atop \nu\neq i}^{d} (T-\lambda_\nu \mathrm{Id}_V).
\] 
In the above product notice  that $T- \lambda_i \mathrm{Id}_V$, $T- \lambda_j \mathrm{Id}_V$ are commuting endomorphisms.  
Since the minimal polynomial of $T$ has degree $d$ all $R$-modules 
$\mathrm{Ker}\Pi_i$ are strictly less than $V$. Moreover there is an element $E$ such that $E \not \in \mathrm{Ker}(\Pi_i)$ for all $1\leq i \leq d$. 
Consider a relation 
\begin{equation}
\label{eq:linearDep}
\sum_{\mu=0}^{d} \gamma_\mu \prod_{\nu=0}^{\mu}(T-\lambda_\mu \mathrm{Id}_V)E,
\end{equation}
where $\prod_{\nu=0}^0 (T- \lambda_\nu \mathrm{Id}_V)E=E$. We fist apply the operator $\prod_{\nu=2}^{d} (T- \lambda_\nu \mathrm{Id}_V)$ to eq. (\ref{eq:linearDep}) and we obtain
\[
0=\gamma_0 \Pi_1 E, 
\]
and by the selection of $E$ we have that $a_0=0$. 
We now apply $\prod_{\nu=3}^{d}(T- \lambda_\nu \mathrm{Id}_V)$ to eq. (\ref{eq:linearDep}). We obtain that
\[
0= \gamma_1 \prod_{\nu=3}^{d}(T- \lambda_\nu \mathrm{Id}_V) (T- \lambda_1 \mathrm{Id}_V) = \gamma_1 \Pi_2 E,
\]
and by the selection of $E$ we have that $\gamma_1=0$. We now apply $\prod_{\nu=4}^{d}(T- \lambda_\nu \mathrm{Id}_V)$ to eq. (\ref{eq:linearDep}) and we obtain
\[
0 = \gamma_2\prod_{\nu=4}^{d}(T- \lambda_\nu \mathrm{Id}_V) (T-\lambda_2 \mathrm{Id}_V)(T-\lambda_1\mathrm{Id}_V) E= \gamma_2 \Pi_3E
\]
and by the selection of $E$ we obtain $\gamma_3=0$. Continuing this way we finally arrive at $\gamma_0=\gamma_1=\cdots=\gamma_{d-1}=0$.
\end{proof}

\begin{lemma}
\label{lemma:Jordan}
Let $V$ be a free $R$-module of rank $R$ acted on by an automorphism
$T:V \rightarrow V$ of order $p^h$. Assume that the minimal polynomial of $T$ is of degree $d$ and has roots $\lambda_1,\ldots,\lambda_d$. 
Then 
$T$ can be written as a matrix with respect to the basis as
follows:
\begin{equation}
\label{eq-matrix-lift}
\begin{pmatrix}
\lambda_1 & 0 & \cdots & \cdots  &0  \\
a_1 &  \lambda_2 &  \ddots &      &\vdots \\
0 &  a_2 & \lambda_3 & \ddots  & \vdots \\
\vdots & \ddots & \ddots & \ddots  & 0 \\
0 &\cdots & 0 & a_{d-1} & \lambda_d
\end{pmatrix}
\end{equation}
\end{lemma}
\begin{proof}
By lemma \ref{lemma:lin-inde} the elements 
\[
E, (T- \lambda_1\mathrm{Id}_V)E, (T- \lambda_2\mathrm{Id}_V)(T- \lambda_1\mathrm{Id}_V)E, \ldots, (T- \lambda_{d-1}\mathrm{Id}_V)\cdots (T- \lambda_1 \mathrm{Id}_V)E\]
form a free submodule of $V$ of rank $d$. 
The theory of submodules of principal ideal domains, there is a basis $E_1,E_2,\ldots, E_{d}$ of the free module $V$ such that 
\begin{align}
\label{eq-E1E2bas}
E_1&=E,
\\
a_1 E_2 &= (T- \lambda_1 \mathrm{Id}_V)E_1,  \nonumber
\\
a_2E_3 &=  (T- \lambda_2\mathrm{Id}_V) E_2, \nonumber\\
&\ldots \nonumber
\\
a_{s-1}E_d&=  (T- \lambda_{d-1}\mathrm{Id}_V) E_{d-1}. \nonumber
\end{align}
Let us consider the module $V_1=\langle E_1,\ldots, E_d \rangle \subset V$. 
By construction, the map $T$ restricts to an automorphism $V_1\rightarrow V_1$ with respect to the basis $E_1,\ldots,E_d$
has the desired form.
We then consider the free module $V/V_1$ and we repeat the procedure for the minimal polynomial of $T$, which again acts on $V/V_1$. The desired result follows. 
\end{proof}
\begin{remark}
The element $T$ as defined in eq. (\ref{eq-matrix-lift}) has order equal to the higher order of the eigenvalues $\lambda_1,\ldots, \lambda_d$ involved. Indeed, since we have assumed that the eigenvalues are different the matrix is diagonalizable in $\mathrm{Quot(R)}$ and has order equal to the maximal order of the eigenvalues involved. In particular it has order $q$ if there is at least one $\lambda_i$ that is a primitive $q$-root of unity. The statement about the order of $T$ is not necessarily true if some of the eigenvalues are the same. For instance the matrix 
$\begin{pmatrix}
1 & 0 \\ 1 & 1  
\end{pmatrix}  
$ has infinite order over a field of characteristic zero. 
\end{remark}
\begin{remark}
The number of indecomposable $R[T]$-summands of $V$ is given by $\#\{i: a_i=0\}+1$.

A lift of a sum of indecomposable $kC_q$-modules $J_{\kappa_1}\oplus \cdots \oplus J_{\kappa_n}$ can form an indecomposable $RC_q$-module. For example the indecomposable module where the generator $T$ of $C_q$ has the form 
\[
T=
\begin{pmatrix}
\lambda_1 & 0 & \cdots & \cdots  &0  \\
a_1 &  \lambda_2 &  \ddots &      &\vdots \\
0 &  a_2 & \lambda_3 & \ddots  & \vdots \\
\vdots & \ddots & \ddots & \ddots  & 0 \\
0 &\cdots & 0 & a_{s-1} & \lambda_d
\end{pmatrix}
\]
where $a_1=\cdots=a_{\kappa_1-1}=1$, $a_{\kappa_1} \in \mathfrak{m}_R$, 
$a_{\kappa_1+1},\ldots, a_{\kappa_2+\kappa_1-1}=1$, $a_{\kappa_2+\kappa_1} \in \mathfrak{m}_R$
, etc
reduces to a decomposable direct sum of Jordan normal forms of sizes 
$J_{\kappa_1},J_{\kappa_2- \kappa_1},\cdots$. 
\end{remark}

\begin{remark}
It is an interesting question to classify these matrices up to conjugation with a matrix in $\mathrm{GL}_{d}(R)$. It seems that the valuation of elements $a_i$ should also play a role.  
\end{remark}


\begin{definition}
Let $h_i(x_1,\ldots,x_j)$ be the complete symmetric polynomial of degree $i$ in the variables $x_1,\ldots,x_j$. For instance
\[
h_3(x_1,x_2,x_3)=x_1^3+x_1^2 x_2+x_1^2 x_3+x_1 x_2^2+x_1
   x_2 x_3+x_1 x_3^2+x_2^3+x_2^2 x_3+x_2
   x_3^2+x_3^3.
\]
Set 
\begin{align*}
L(\kappa,j,\nu)&=h_\kappa(\lambda_j,\lambda_{j+1},\ldots,\lambda_{j+\nu})\\
A(i,j)&= 
\begin{cases}
a_i a_{i+1} \cdots a_{i+j} & \text{ if } j \geq 0 \\
0 & \text{ if } j<0
\end{cases}
\end{align*}
\end{definition}
\begin{lemma}
The matrix $T^\alpha=(t_{ij}^{(\alpha)})$ is given by the following formula:
\[
t_{ij}^{(\alpha)}=
\begin{cases}
\lambda_i^\alpha  & \text{ if } i=j\\
A(j,i-j-1) \cdot  L(\alpha-(i-j),j,i-j) & \text{ if } j <  i \\
0  & \text{ if } j > i
\end{cases}
\]
\end{lemma}
\begin{proof}
For $j\geq i$ the proof is trivial. When $j< i$ and $\alpha=1$ it is immediate, since $L(x,\cdot,\cdot)\equiv0$, for every $x\leq0$. Assume this holds for $\alpha=n$. If $\alpha=n+1$, 
\begin{align*}
  t_{ij}^{(n+1)} &= t_{ij}^{(n)}t_{ij} = \sum_{k=1}^{r}t_{ik}^{(\alpha)}t_{kj}=\lambda_j t^{(\alpha)}_{ij}+a_j t^{(\alpha)}_{ij+1} = \lambda_j A(j,i-j-1)L(\alpha-(i-j),j,i-j) +\\
                 &+ a_jA(j+1,i-j-2)L(\alpha-(i-j-1),j+1,i-j-1) = \\
                 &= A(j,i-j-1)\left(\lambda_j h_{\alpha-(i-j)}(\lambda_j, \dots, \lambda_{j})+ h_{\alpha-(i-j)+1}(\lambda_{j+1},\dots, \lambda_i) \right) = \\
                 &= A(j,i-j-1)h_{\alpha-(i-j)+1}(\lambda_{j},\dots, \lambda_i) =\\
                 &= A(j,i-j-1)L(\alpha-(i-j)+1, i, i-j).
\end{align*}
\end{proof}

\begin{remark}
The space of homogeneous polynomials of degree $k$ in $n$-variables has dimension $\binom{n-1+c}{n-1}$. Since all $q$-roots of unity are reduced to $1$ modulo $\mathfrak{m}_R$ the quantity 
$L(\alpha-(i-j),j,i-j)$ is reduced to $n=(i-j)+1$, $c=\alpha-(i-j)$
\[
\binom{n-1+c}{n-1}=
\binom{\alpha}{i-j}. 
\]
This equation is compatible with the computation of $\tau^\alpha$ given in eq. (\ref{eq:a-power-t-spec}). 
\end{remark}
\begin{lemma}
\label{lemma:eigenvector-s-sel}
There is an eigenvector $E$ of the generator $\sigma$ of the cyclic group $C_m$ which is not an element in
 $\displaystyle\bigcup_{i=1}^s \mathrm{Ker}(\Pi_i\otimes K)$.
\end{lemma}
\begin{proof}
The eigenvectors $E_1,\ldots,E_d$ of $\sigma$ form a basis of the space $V\otimes K$. By multiplying by certain elements in $R$, if necessary, we can assume that all $E_i$ are in $V$ and their reductions $E_i\otimes R/\mathfrak{m}_R$, $1\leq i \leq d$ give rise to a basis of eigenvectors of a generator of the cyclic group $C_m$ acting on $V\otimes R/\mathfrak{m}_R$. If every eigenvector $E_i$ is an element of some $\mathrm{Ker}(\Pi_\nu)$ for $1\leq i \leq d$, then their reductions will be elements in $\mathrm{Ker}(T-1)^{d-1}$, a contradiction since the later kernel has dimension $<d$. 
\end{proof}

\begin{lemma}
\label{lemma:recursive}
Let $V$ be a free $C_q \rtimes C_m$-module, which is indecomposable as a $C_q$-module. 
Consider the basis given in lemma \ref{lemma:Jordan}. Then the value of $\sigma(E_1)$ determines $\sigma(E_i)$ for $2\leq i \leq d$. 
\end{lemma}
\begin{proof}
Let $\sigma$ be a generator of the cyclic group $C_m$.  We will use the notation of lemma \ref{lemma:lin-inde}.  
We use lemma \ref{lemma:eigenvector-s-sel} in order to select a suitable eigenvector of $E_1$ of $\sigma$
and then form the basis $E_1,E_2,\ldots,E_d$ as given in eq. (\ref{eq-E1E2bas}). 
We can compute the action of $\sigma$ on all basis elements $E_i$ by
\begin{equation}
\label{eq:ind-s}
\sigma(a_{i-1}E_i)= \sigma(T- \lambda_{i-1}  \mathrm{Id}_V)E_{i-1}=
(T^a - \lambda_{i-1} \mathrm{Id}_V)\sigma(E_{i-1}). 
\end{equation}
This means that one can define recursively the action of $\sigma$ 
on all elements $E_i$. 
Indeed, 
assume that 
\[
\sigma(E_{i-1}) =\sum_{\nu=1}^{d} \gamma_{\nu,i-1} E_{\nu}. 
\]
We now have
\begin{align*}
(T^a-\lambda_{i-1} \mathrm{Id}_V) E_{\nu} &= 
\sum_{\mu=1}^{d} t_{\mu,\nu}^{(\alpha)} E_\mu - \lambda_{i-1} E_\nu
\\
&= (\lambda_{\nu}^\alpha - \lambda_{i-1}) E_{\nu}+ \sum_{\mu=\nu+1}^{d} t_{\mu,\nu}^{(\alpha)} E_\mu
\end{align*}
We combine all the above to 
\begin{align}
\nonumber
a_{i-1} \sigma(E_i) &= \sum_{\nu=1}^d \gamma_{\nu,i-1}(\lambda_\nu^\alpha-\lambda_{i-1})E_\nu 
+
\sum_{\nu=1}^d \gamma_{\nu,i-1}\sum_{\mu=\nu+1}^d  t_{\mu,\nu}^{(\alpha)} E_\mu
\\
\label{eq:allcombined}
&=\sum_{\nu=1}^d \tilde{\gamma}_{\nu,i} E_\nu,
\end{align}
for a selection of elements $\gamma_{\nu,i}\in R$, which can be explicitly computed by collecting the coefficients of the basis elements $E_1,\ldots,E_d$. 

Observe that the quantity on the right hand side of eq. (\ref{eq:allcombined}) must be divisible by $a_{i-1}$. 
Indeed, 
 let $v$ be the valuation of the local principal ideal domain $R$. 
Set 
\[
e_0= \min_{1\leq \nu \leq d}\{
v(\tilde{\gamma}_{\nu,i})
\}.
\]
If $e_0< v(a_{i-1})$ then we 
 divide eq. (\ref{eq:allcombined}) by $\pi^{e_0}$ where $\pi$ is the local uniformizer of $R$, that is $\mathfrak{m}_R=\pi R$. We then consider the divided equation modulo $\mathfrak{m}_R$ to obtain a linear dependence relation among the elements 
$E_i\otimes k$, which is a contradiction. 
Therefore $e_0 \geq v(a_{i-1})$ and we obtain an 
equation 
\[
\sigma(E_i) = \sum_{\nu=1}^d
\frac{
\tilde{\gamma}_{\nu,i}
}{
  a_{i-1}
} E_\nu
=
\sum_{\nu=1}^d 
\gamma_{\nu,i}
 E_\nu.
\]
\end{proof}
For example $\sigma(E_1)=\zeta_m^{\epsilon}E_1$. We compute that 
\[
a_1 \sigma(E_2)= (T^\alpha-\lambda_1\mathrm{Id}) \sigma(E_1)
\]
and
\begin{align*}
\sigma(E_2) &=\frac{(\lambda_1^\alpha -\lambda_1)}{a_1} \zeta_\mu^{\epsilon} E_1 +
\zeta_m^{\epsilon} \sum_{\mu=2}^d \frac{t_{\mu,1}^{(\alpha)}}{a_1}E_\mu
\\
&=\frac{(\lambda_1^\alpha -\lambda_1)}{a_1} \zeta_\mu^{\epsilon} E_1 +
\zeta_m^{\epsilon} \sum_{\mu=2}^d 
\frac{
A(1,\mu-2)
L(\alpha-(\mu-1),1,\mu-1)
}{a_1}E_\mu
\\
&=\frac{(\lambda_1^\alpha -\lambda_1)}{a_1} \zeta_\mu^{\epsilon} E_1 +
\zeta_m^{\epsilon} \sum_{\mu=2}^d
\frac{
a_1a_2\cdots a_{\mu-1} 
h_{\alpha-(\mu-1)}(\lambda_1,\lambda_2,\ldots,\lambda_{\mu})
}{a_1}E_\mu.
\end{align*}

\begin{proposition}
\label{prop:fing-G-uniqe}
Assume that no element $a_1,\ldots,a_{d-1}$ given in eq. (\ref{eq-matrix-lift}) is zero. 
Given $\alpha\in \mathbb{N}, \alpha\geq 1$ and an element $E_1$, which is not an element in $\bigcup_{i=1}^{d} \mathrm{Ker}(\Pi_i\otimes K)$, if there is a  matrix $\Gamma=(\gamma_{ij})$, such that $\Gamma T \Gamma^{-1}=T^\alpha$ and $\Gamma E_1= \zeta_m^\epsilon E_1$, then this matrix $\Gamma$ is unique.   
\end{proposition}
\begin{proof}
We will use the idea leading to equation (\ref{eq:ind-s}) replacing $\sigma$ with $\Gamma$. We will compute recursively and uniquely the entries $\gamma_{\mu,i}$, arriving at the explicit formula of eq. (\ref{eq:fin-sum4}).

Observe that trivially $\gamma_{\nu,1}=0$ for all $\nu<1$ since we only allow $1\leq \nu \leq d$. We compute
\begin{align}
\label{next-gamma-i}
\tilde{\gamma}_{\mu,i} &=
\gamma_{\mu,i-1}(\lambda_\mu^\alpha- \lambda_{i-1}) + \sum_{\nu=1}^{\mu-1} \gamma_{\nu,i-1} t_{\mu,\nu}^{(\alpha)} \\
&= 
\gamma_{\mu,i-1}(\lambda_\mu^\alpha- \lambda_{i-1}) + \sum_{\nu=1}^{\mu-1} \gamma_{\nu,i-1} 
A(\nu, \mu-\nu-1) L\big(\alpha-(\mu-\nu),\nu,\mu-\nu)
\nonumber
\\
&= 
\gamma_{\mu,i-1}(\lambda_\mu^\alpha- \lambda_{i-1}) + \sum_{\nu=1}^{\mu-1} \gamma_{\nu,i-1} 
a_{\nu}a_{\nu+1}\cdots a_{\mu-1}
h_{\alpha- \mu+\nu}(\lambda_\nu, \lambda_{\nu+1},\ldots,\lambda_\mu)
\nonumber
\end{align}

Define 
\begin{align*}
[\lambda_m^\alpha- \lambda_x]_{i}^j &=\prod_{x=i}^j (\lambda_\mu^ \alpha- \lambda_x)
\\
[a]_i^j &= \prod_{x=i}^j a_x
\end{align*}
for $i \leq j$. If $i>j$ then both of the above quantities are defined to be equal to $1$. 

Observe that for  $\mu=1$ eq. (\ref{next-gamma-i}) becomes  
\begin{equation}
\label{eq:gamma1}
  \gamma_{1,i}=\frac{1}{a_{i-1}}\gamma_{1,i-1}(\lambda_1^\alpha - \lambda_{i-1})
\end{equation}
and we arrive at (assuming that $\Gamma(E_1)=\zeta_m^{\epsilon} E_1$)
\begin{equation}
\label{eq:1s}
  \gamma_{1,i} =  \frac{\zeta_m^{\epsilon}}{a_1a_2\cdots a_{i-1}}\prod_{x=1}^{i-1} (\lambda_1^\alpha - \lambda_x)=
  \frac{\zeta_m^\epsilon}{a_1a_2\cdots a_{i-1}} [\lambda_1^\alpha- \lambda_x]_1^{i-1}. 
\end{equation}
For $\mu \geq 2$ we have $\gamma_{\mu,1}=0$, since by assumption $TE_1=\zeta_m^{\epsilon} E_1$.  Therefore eq. (\ref{next-gamma-i}) gives us
\begin{align}
\nonumber
  \gamma_{\mu,i} &=
  \sum_{\kappa_1=0}^{i-2}  
  \frac{
  [\lambda_\mu^\alpha - \lambda_x]_{i-\kappa_1}^{i-1}
  }
  {
  [a]_{i-1-\kappa_1}^{i-1}
  }
\sum_{\mu_2=1}^{\mu-1} \gamma_{\mu_2,i-1- \kappa_1} [a]_{\mu_2}^{\mu-1} h_{\alpha- \mu+\mu_2}(\lambda_{\mu_2},\ldots, \lambda_\mu)
\\
&=
\label{eq:ind-11}
\sum_{\mu_2=1}^{\mu-1}
  [a]_{\mu_2}^{\mu-1} h_{\alpha- \mu+\mu_2}(\lambda_{\mu_2},\ldots, \lambda_\mu)
  \sum_{\kappa_1=0}^{i-2}
  \frac{
  [\lambda_\mu^\alpha - \lambda_x]_{i-\kappa_1}^{i-1}
  }
  {
  [a]_{i-1-\kappa_1}^{i-1}
  }
  \gamma_{\mu_2,i-1-\kappa_1}.
\end{align}
We will now prove eq. (\ref{eq:ind-11}) by induction on $i$. 
For $i=2,\mu \geq 2$ we have 
\begin{align*}
  \gamma_{\mu,2} &= \frac{1}{a_1} \gamma_{\mu,1}(\lambda_\mu^\alpha- \lambda_1) + 
  \frac{1}{a_1} \sum_{\mu_2=1}^{\mu-1} \gamma_{\mu_2,1} [a]_{\mu_2}^{\mu-1} h_{\alpha-\mu+\mu_2}(\lambda_{\mu_2},\ldots,\lambda_\mu)
  \\
  &= \frac{1}{a_1} [a]_1^{\mu-1} h_{\alpha-\mu+1}(\lambda_1,\ldots,\lambda_\mu) \gamma_{1,1}.
\end{align*}
Assume now that eq. (\ref{eq:ind-11})  holds for computing  $\gamma_{\mu,i-1}$. We will treat the $\gamma_{\mu,i}$ case. We have 
\begin{align*}
\gamma_{\mu,i} &= \frac{(\lambda_\mu^\alpha- \lambda_{i-1})}{a_{i-1}} \gamma_{\mu,i-1} + 
\frac{1}{a_{i-1}}
\sum_{\mu_2=1}^{\mu-1} \gamma_{\mu_2,i-1} [a]_{\mu_2}^{\mu-1} h_{\alpha-\mu+\mu_2}(\lambda_{\mu_2},\ldots,\lambda_\mu)
\\
&=
\frac{(\lambda_\mu^\alpha- \lambda_{i-1})}{a_{i-1}}
\sum_{\mu_2=1}^{\mu-1}
  [a]_{\mu_2}^{\mu-1} h_{\alpha- \mu+\mu_2}(\lambda_{\mu_2},\ldots, \lambda_\mu)
  \sum_{\kappa_1=0}^{i-3}
  \frac{
  [\lambda_\mu^\alpha - \lambda_x]_{i-1-\kappa_1}^{i-2}
  }
  {
  [a]_{i-2-\kappa_1}^{i-2}
  }
  \gamma_{\mu_2,i-2-\kappa_1}
  \\
  &+\frac{1}{a_{i-1}}
\sum_{\mu_2=1}^{\mu-1} \gamma_{\mu_2,i-1} [a]_{\mu_2}^{\mu-1} h_{\alpha-\mu+\mu_2}(\lambda_{\mu_2},\ldots,\lambda_\mu)
\\
&=
\sum_{\mu_2=1}^{\mu-1}
  [a]_{\mu_2}^{\mu-1} h_{\alpha- \mu+\mu_2}(\lambda_{\mu_2},\ldots, \lambda_\mu)
  \sum_{\kappa_1=0}^{i-3}
  \frac{
  [\lambda_\mu^\alpha - \lambda_x]_{i-1-\kappa_1}^{i-1}
  }
  {
  [a]_{i-2-\kappa_1}^{i-1}
  }
  \gamma_{\mu_2,i-2-\kappa_1}
  \\
&+
  \frac{1}{a_{i-1}}
\sum_{\mu_2=1}^{\mu-1} \gamma_{\mu_2,i-1} [a]_{\mu_2}^{\mu-1} h_{\alpha-\mu+\mu_2}(\lambda_{\mu_2},\ldots,\lambda_\mu)
\\
&=
\sum_{\mu_2=1}^{\mu-1}
  [a]_{\mu_2}^{\mu-1} h_{\alpha- \mu+\mu_2}(\lambda_{\mu_2},\ldots, \lambda_\mu)
  \sum_{\kappa_1=1}^{i-2}
  \frac{
  [\lambda_\mu^\alpha - \lambda_x]_{i-\kappa_1}^{i-1}
  }
  {
  [a]_{i-1-\kappa_1}^{i-1}
  }
  \gamma_{\mu_2,i-1-\kappa_1}
  \\
&+
\sum_{\mu_2=1}^{\mu-1} [a]_{\mu_2}^{\mu-1} h_{\alpha-\mu+\mu_2}(\lambda_{\mu_2},\ldots,\lambda_\mu)
  \frac{1}{a_{i-1}}
  \gamma_{\mu_2,i-1}
\\
&=
\sum_{\mu_2=1}^{\mu-1}
  [a]_{\mu_2}^{\mu-1} h_{\alpha- \mu+\mu_2}(\lambda_{\mu_2},\ldots, \lambda_\mu)
  \sum_{\kappa_1=0}^{i-2}
  \frac{
  [\lambda_\mu^\alpha - \lambda_x]_{i-\kappa_1}^{i-1}
  }
  {
  [a]_{i-1-\kappa_1}^{i-1}
  }
  \gamma_{\mu_2,i-1-\kappa_1}
\end{align*}
and equation (\ref{eq:ind-11}) is now proved. 

We proceed recursively applying eq. (\ref{eq:ind-11}) to each of the summands $\gamma_{\mu_2,i-1-\kappa_1}$ if 
$\mu_2 >1$ and $i-1-\kappa_1 >1$. If $\mu_2=1$, then $\gamma_{\mu_2,i-1-\kappa_1}$ is computed by eq. (\ref{eq:gamma1}) and if $\mu_2>1$ and $i-1-\kappa_1 \leq 1$ then $\gamma_{\mu_2,i-1-\kappa_1}=0$. 
We can classify all iterations needed by the set $\Sigma_\mu$ of sequences $(\mu_s,\mu_{s-1},\ldots, \mu_3,\mu_2)$ such that 
\begin{equation}
\label{eq:mu-iter}
1=\mu_s < \mu_{s-1} < \cdots < \mu_3 < \mu_2 < \mu=\mu_1.
\end{equation}
For example for $\mu=5$ the set of such sequences is given by 
\[
  \Sigma_\mu=\{ (1), (1,2), (1,3), (1,2,3), (1,4), (1,2,4), (1,3,4), (1,2,3,4)\}
\]
corresponding to the  tree of iterations given in figure \ref{fig:Iter-Tree}. 
The length of the sequence $(\mu_s,\mu_{s-1},\ldots,\mu_2)$ is  given in eq. (\ref{eq:mu-iter}) is $s-1$. In each iteration the $i$ changes to $i-1-k$ thus we have the following sequence of indices 
\[
i_1=i\rightarrow i_2=i-1-\kappa_1 \rightarrow i_3=i-2-(\kappa_1+\kappa_2) \rightarrow 
\cdots
\rightarrow 
i_{s} = i-(s-1) - (\kappa_1+\cdots+ \kappa_{s-1})
\] 
\begin{figure}
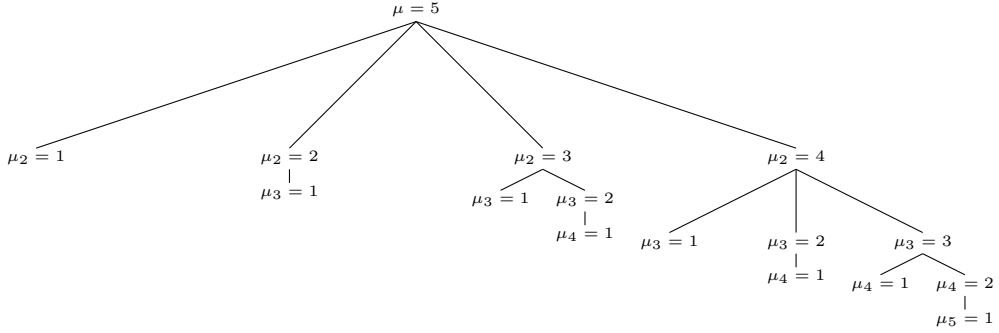

{ \tiny
\Tree[.$\mu=5$ [.$\mu_2=1$ ]
          [.$\mu_2=2$ [.$\mu_3=1$ ]
                ] 
          [.$\mu_2=3$ [.$\mu_3=1$ ] [.$\mu_3=2$ [.$\mu_4=1$ ]] ]
          [.$\mu_2=4$ [.$\mu_3=1$ ]  [.$\mu_3=2$ [.$\mu_4=1$ ] ] [.$\mu_3=3$ [.$\mu_4=1$ ] [.$\mu_4=2$ 
           [.$\mu_5=1$ ]
          ]]]
                  ]
                  }
\caption{ \label{fig:Iter-Tree}  Iteration tree for $\mu=5$}
\end{figure}
For the sequence  $i_1,i_2,\ldots,$ we might have $i_t=1$ for $t<s-1$. But in this case, we will arrive at the element $\gamma_{\mu_{t+1},i_t}=\gamma_{\mu_t,1}=0$ since $\mu_t>1$. 
This means that we will have to consider only selections $\kappa_1,\ldots,\kappa_{s-1}$ such that $i_{s-1} \geq 1$. Therefore we arrive at the following
 expression for $\mu\geq 2$
\begin{align}
\nonumber
\gamma_{\mu,i} &=
\sum_{(\mu_s,\ldots,\mu_2)\in \Sigma_\mu}
[a]_{\mu_2}^{\mu-1}[a]_{\mu_3}^{\mu_2-1}\cdots [a]_{\mu_s}^{\mu_{s-1}-1}
\prod_{\nu=2}^s h_{\alpha-\mu_{\nu-1}+ \mu_\nu}(\lambda_{\mu_\nu},\ldots,\lambda_{\mu_{\nu-1}})
\\ \nonumber
&\cdot
\sum_{i=i_1> i_2 > \cdots > i_{s} \geq 1} 
\prod_{\nu=1}^{s-1}
\frac{
  [\lambda_{\mu_\nu}^\alpha- \lambda_x]^{i_\nu-1}_{i_{\nu+1}+1}
}
{
  [a]^{i_\nu-1}_{i_{\nu+1}}
}
  \cdot
  \gamma_{1,i_s}.
  \\
&=
\nonumber
\sum_{(\mu_s,\ldots,\mu_2)\in \Sigma_\mu}
\prod_{\nu=2}^s h_{\alpha-\mu_{\nu-1}+ \mu_\nu}(\lambda_{\mu_\nu},\ldots,\lambda_{\mu_{\nu-1}})
\\
\nonumber
&\cdot
\sum_{i=i_1> i_2 > \cdots > i_{s} \geq 1} 
\frac{
  [a]_1^{\mu-1}
}
{
  [a]^{i-1}_{i_{s}}
}
\prod_{\nu=1}^{s-1}
  [\lambda_{\mu_\nu}^\alpha- \lambda_x]^{i_\nu-1}_{i_{\nu+1}+1}
  \frac{ \zeta_m^\epsilon [\lambda_1^\alpha- \lambda_x]_1^{i_{s}-1}}{[a]_1^{i_{s}-1}}
\\
\label{eq:fin-sum4}
&=
\!\!\!\!\!\!\!\!\!\!\!
\sum_{(\mu_s,\ldots,\mu_2)\in \Sigma_\mu}
\prod_{\nu=2}^s
\!\! h_{\alpha-\mu_{\nu-1}+ \mu_\nu}(\lambda_{\mu_\nu},\ldots,\lambda_{\mu_{\nu-1}})
\frac{
  [a]_1^{\mu-1}
}
{
  [a]^{i-1}_{1}
}
\zeta_m^\epsilon
\!\!\!
\sum_{i=i_1> i_2 > \cdots > i_{s} \geq 1} 
\prod_{\nu=1}^{s}
  [\lambda_{\mu_\nu}^\alpha- \lambda_x]^{i_\nu-1}_{i_{\nu+1}+1}
\end{align}
where  $i_{s+1}+1=1$ that is $i_{s+1}=0$.
\end{proof}

We will now prove that the matrix $\Gamma$ of lemma \ref{prop:fing-G-uniqe} exists by cheking 
that $\Gamma T = T^\alpha \Gamma$. 
Set $(a_{\mu,i})=\Gamma T$, $(b_{\mu,i})= T^\alpha \Gamma$. For $i<d$ we have 
\begin{align*}
  a_{\mu,i} &=\sum_{\nu=1}^{d} \gamma_{\mu,\nu} t_{\nu,i} = \gamma_{\mu,i} t_{ii} + \gamma_{\mu,i+1} t_{i+1,i}
  \\
  &= \gamma_{\mu,i} \lambda_i + \gamma_{\mu,i}(\lambda_\mu^\alpha- \lambda_i) + \sum_{\nu=1}^{\mu-1} \gamma_{\nu,i} t_{\mu,\nu}^{(\alpha)} \\
  &=  \gamma_{\mu,i} \lambda_\mu^\alpha + \sum_{\nu=1}^{\mu-1} \gamma_{\nu,i} t_{\mu,\nu}^{(\alpha)}
= 
  \sum_{\nu=1}^{\mu}  t_{\mu,\nu}^{(\alpha)} \gamma_{\nu,i} = b_{\mu,i}. 
\end{align*}
For $i=d$ we have:
\begin{align*}
  a_{\mu,d} &=\sum_{\nu=1}^d \gamma_{\mu,\nu} t_{\nu,d} = \gamma_{\mu,d} t_{d,d} =\gamma_{\mu,d} \lambda_d
\end{align*}
while 
\begin{align*}
b_{\mu,d}&=\sum_{\nu=1}^d t_{\mu,\nu}^{(\alpha)} \gamma_{\nu,d}= \sum_{\nu=1}^{\mu-1} t_{\mu,\nu}^{(\alpha)} \gamma_{\nu,d} + \lambda_\mu^{\alpha} \gamma_{\mu,d}
\end{align*}
This gives us the relation
\begin{equation}
\label{eq:recf}
(\lambda_d- \lambda_\mu^a) \gamma_{\mu,d}= \sum_{\nu=1}^{\mu-1} t_{\mu,\nu}^{(\alpha)} \gamma_{\nu,d}
\end{equation}
For $\mu=1$ using eq. (\ref{eq:1s}) we have
\[
  \gamma_{1,d} \lambda_d = \gamma_{1,d} \lambda_1^\alpha \Rightarrow 
[\lambda_1^\alpha - \lambda_x]_1^{d} = 0.
\]
This relation is satisfied if $\lambda_1^\alpha$ is one of $\{\lambda_1,\ldots,\lambda_d\}$. 
Without loss of generality we assume that 
\begin{equation}
\label{eq:lambda-act}
  \lambda_i^{(a)}=
  \begin{cases}
  \lambda_{i+1} & \text{ if } m \nmid i \\
  \lambda_{i-m+1} & \text{ if } m \mid i  
  \end{cases}
\end{equation}
We have the following conditions:
\begin{align*}
\mu=2 & &(\lambda_d- \lambda_2^\alpha) \gamma_{2,d} & = t_{2,1}^{(\alpha)} \gamma_{1,d}
\\
\mu=3 & & (\lambda_d- \lambda_3^\alpha) \gamma_{3,d} & = 
  t_{3,1}^{(\alpha)} \gamma_{1,d} +
  t_{3,2}^{(\alpha)} \gamma_{2,d} 
  \\
\mu=4 && 
 (\lambda_d- \lambda_4^\alpha) \gamma_{4,d} &= 
  t_{4,1}^{(\alpha)} \gamma_{1,d} +
  t_{4,2}^{(\alpha)} \gamma_{2,d} +
  t_{4,3}^{(\alpha)} \gamma_{3,d}
  \\
\vdots && \vdots 
\\
\mu=d-1 && 
(\lambda_d- \lambda_{d-1}^\alpha) \gamma_{d-1,d} &= 
  t_{d-1,1}^{(\alpha)} \gamma_{1,d} +
  t_{d-1,2}^{(\alpha)} \gamma_{2,d} + \cdots +
  t_{d-1,d-2}^{(\alpha)} \gamma_{d-1,d}
\end{align*}
All these equations are true provided that $\gamma_{1,d},\ldots,\gamma_{d-2,d}=0$. 
Finally, for $\mu=d$, we have
\begin{equation}\label{finalEquation}
  (\lambda_d - \lambda_d^{\alpha}) \gamma_{d,d} =  \sum_{\nu=1}^{d-1} t_{d,\nu}^{(\alpha)} \gamma_{\nu,d}
\end{equation}
which is true provided that $(\lambda_d - \lambda^\alpha_d) \gamma_{d,d}= t_{d,d-1}^{(a)} \gamma_{d-1,d}.$

\begin{lemma}
\label{lemma:sum-prod}
For $n\geq 2$ the vertical sum $S_n$ of the products of every line  of the following array 
\[
\begin{array}{c|cccccc}
y\\
\hline 
1 &  1  & (x_1-x_2) & (x_1-x_3) & \cdots & \cdots & (x_1-x_{n}) \\
2 &  (z-x_1) & 1  & (x_1-x_3) & \cdots & \cdots & (x_1-x_{n}) \\
3 &  (z-x_1) & (z-x_2) & 1  & \ddots & \ddots & \vdots \\
\vdots &\vdots  & \ddots & \ddots & \ddots &  & \vdots 
\\
\vdots &\vdots  & &  & \ddots &  &   \vdots 
\\
n-1 & (z-x_1) & (z-x_2) & \cdots  &(z-x_{n-2}) & 1 & (x_1-x_{n}) 
\\
n &  (z-x_1) & (z-x_2) & \cdots &  (z-x_{n-2})  & (z-x_{n-1}) &  1 
\end{array}
\]
is given by
\[S_n=
\sum_{y=1}^n  \prod_{\nu=y+1}^{n} (x_1-x_\nu) \prod_{\mu=1}^{y-1} (z- x_\mu) 
 = (z-x_2) \cdots (z - x_{n}). 
\] 
In particular when $z=x_{n}$ the sum is zero. 
\end{lemma}
\begin{proof}
We will prove the lemma by induction. 
For $n=2$ we have $S_2=(x_1-x_2)+ (z-x_1)=z-x_2$. Assume that the equality holds for $n$. The sum $S_{n+1}$ corresponds to the array:
\[
\begin{array}{c|cccccc}
y\\
\hline 
1 &  1  & (x_1-x_2) & (x_1-x_3) & \cdots & (x_1-x_n) & (x_1-x_{n+1})  \\
2 &  (z-x_1) & 1  & (x_1-x_3) & \cdots & (x_1- x_n) & (x_1 -x_{n+1})\\
3 &  (z-x_1) & (z-x_2) & 1  & \ddots &  \vdots &  \vdots\\
\vdots &\vdots  &  & \ddots & \ddots   & \vdots & \vdots 
\\
n-1 & (z-x_1) & \cdots & (z-x_{n-2}) & 1 & (x_1-x_n) & (x_1-x_{n+1})
\\
n &  (z-x_1) & (z-x_2) & \cdots & (z-x_{n-1}) &  1 & (x_1-x_{n+1}) \\
n+1 & (z-x_1) & (z-x_2) & \cdots & (z-x_{n-1}) & (z-x_n)&  1  
\end{array}
\]
We have by definition $S_{n+1}=S_n (x_1-x_{n+1}) + (z-x_1)(z-x_2)\cdots (z-x_{n})$, which by induction gives 
\begin{align*}
S_{n+1} &=  (z -x_2)\cdots(z-x_{n})(x_1-x_{n+1}) + (z-x_1)(z-x_2)\cdots (z-x_{n}) 
\\
&=   (z -x_2)\cdots(z-x_{n})(x_1-x_{n+1}+z-x_1)
\end{align*}
and gives the desired result.
\end{proof}

\begin{lemma}
\label{lemma:ABlL}
Consider $A < l < L <B$. The quantity
\[
 \sum_{l \leq y \leq L} 
 [\lambda_{a} - \lambda_x]_{A}^{y-1}
 \cdot
 [\lambda_{b} - \lambda_x]_{y+1}^{B} 
\]
equals to
\[
[\lambda_{a} - \lambda_x]_{A}^{l-1}
 \cdot
 [\lambda_{b} - \lambda_x]_{L+1}^{B}
  \cdot 
\frac{
  [\lambda_{a}- \lambda_x]_{l}^{L} -  [\lambda_{b} - \lambda_x]_l^L
}
{(\lambda_{a}- \lambda_{b})}
\]
\end{lemma}
\begin{proof}
We write
\[
 \sum_{l \leq y \leq L}
 [\lambda_{a} - \lambda_x]_{A}^{y-1}
 \cdot
 [\lambda_{b} - \lambda_x]_{y+1}^{B} 
\]
\[
 = 
[\lambda_{a} - \lambda_x]_{A}^{l-1}
 \cdot
 [\lambda_{b} - \lambda_x]_{L+1}^{B}
  \cdot 
\sum_{l \leq y \leq L}
[\lambda_{a} - \lambda_x]_{l}^{y-1}
 \cdot
 [\lambda_{b} - \lambda_x]_{y+1}^{L}
 \]
The last sum can be read as the vertical sum $S$ of the products of every line in the following array:
\[
{
\setlength{\arraycolsep}{0.3pt}
\begin{array}{c|ccccccc}
y & 
\\
\hline
l & 1 &
(\lambda_{b}- \lambda_{l+1}) &
(\lambda_{b}- \lambda_{l+2}) &  \cdots & 
(\lambda_{b}- \lambda_{L-1}) 
&  (\lambda_{b}- \lambda_{L})
\\
l+1 & (\lambda_{a} - 
\lambda_{l} ) 
& 1 &
 (\lambda_{b}- \lambda_{l+2}) &  \cdots & 
(\lambda_{b}- \lambda_{{L-1}})
&  (\lambda_{b}- \lambda_{L})
\\
l+2 & 
(\lambda_{a} - \lambda_{l} ) &  (\lambda_{a}- \lambda_{l+1})
& 1 & 
&  \vdots & \vdots 
\\
\vdots  & \vdots & \vdots  & \ddots & \ddots & \vdots & \vdots
\\
L-2& 
(\lambda_{a} - \lambda_{l} ) &  (\lambda_{a}- \lambda_{l+1})
& \cdots & 1 & (\lambda_{b}- \lambda_{{L-1}})
&  (\lambda_{b}- \lambda_{L})
\\
L-1 & 
(\lambda_{a} - \lambda_{l} ) &  (\lambda_{a}- \lambda_{l+1})
& \cdots & 
(\lambda_a- \lambda_{L-2}) & 1
&  (\lambda_{b}- \lambda_{L})
\\
L & 
(\lambda_{a} - \lambda_{l} ) &  (\lambda_{a}- \lambda_{l+1})
& \cdots & (\lambda_a- \lambda_{L-2}) & 
(\lambda_a- \lambda_{{L}-1}) & 
1
\end{array}
}
 \]
If $l=b$, then lemma \ref{lemma:sum-prod} implies that $S=[\lambda_{a}- \lambda_x]_{b+1}^L$. Furthermore, if $L=a$ then $S=0$.

The quantity $S$ cannot be directly computed using  lemma \ref{lemma:sum-prod}, 
if $l\neq b$. We proceed by forming the array: 
\begin{center}
\scalebox{.78}
{
$
\setlength{\arraycolsep}{0.31pt}
\begin{array}{c|ccc|ccccccc}
y & 
\\
\hline
b & 1 & (\lambda_{b}- \lambda_{b+1}) & \cdots 
& 
(\lambda_b - \lambda_{l} ) 
& \cdots  & \cdots & \cdots & \cdots  & (\lambda_{b}-\lambda_L) 
\\
\vdots & & &  & \vdots & & & & & \vdots
\\
l-1 & (\lambda_{a}-\lambda_{b})&\cdots& 1 & (\lambda_b - 
\lambda_{l} ) 
&\cdots  & \cdots & \cdots  & \cdots & (\lambda_{b}-\lambda_L)
\\
\hline
l &
(\lambda_{a}-\lambda_{b})&\cdots&(\lambda_{a}- \lambda_{l-1})&
 1 & 
(\lambda_{b}- \lambda_{l+1}) &
(\lambda_{b}- \lambda_{l+2}) &  \cdots & 
(\lambda_{b}- \lambda_{L-1}) 
&  (\lambda_{b}- \lambda_{L})
\\
l+1 & 
(\lambda_{a}-\lambda_{b})&\cdots&(\lambda_{a}- \lambda_{l-1})&
(\lambda_{a} - 
\lambda_{l} ) 
& 1 &
 (\lambda_{b}- \lambda_{l+2}) &  \cdots & 
(\lambda_{b}- \lambda_{{L-1}})
&  (\lambda_{b}- \lambda_{L})
\\
l+2 & 
(\lambda_{a}-\lambda_{b})&\cdots&(\lambda_{a}- \lambda_{l-1})&
(\lambda_{a} - \lambda_{l} ) &  (\lambda_{a}- \lambda_{l+1})
& 1 & 
 &  \vdots & \vdots 
\\
\vdots 
& & &
 & \vdots & \vdots  & \ddots & \ddots & \vdots & \vdots
\\
L-2& 
(\lambda_{a}-\lambda_{b})&\cdots&(\lambda_{a}- \lambda_{l-1})&
(\lambda_{a} - \lambda_{l} ) &  (\lambda_{a}- \lambda_{l+1})
& \cdots & 1 & (\lambda_{b}- \lambda_{{L-1}})
&  (\lambda_{b}- \lambda_{L})
\\
L-1 & 
(\lambda_{a}-\lambda_{b})&\cdots&(\lambda_{a}- \lambda_{l-1})&
(\lambda_{a} - \lambda_{l} ) &  (\lambda_{a}- \lambda_{l+1})
& \cdots & (\lambda_a- \lambda_{L-2}) & 1
&  (\lambda_{b}- \lambda_{L})
\\
L & 
(\lambda_{a}-\lambda_{b})&\cdots&(\lambda_{a}- \lambda_{l-1})&
(\lambda_{a} - \lambda_{l} ) &  (\lambda_{a}- \lambda_{l+1})
& \cdots & (\lambda_a- \lambda_{L-2}) & 
(\lambda_a- \lambda_{{L}-1}) & 
1
\end{array}
$
}
\end{center}
The value of this array is computed using lemma \ref{lemma:sum-prod} to be equal to 
$
[\lambda_{a} - \lambda_x]_{b+1}^L 
$.
We observe that the sum of the products of the top left array can be computed using lemma \ref{lemma:sum-prod}, while the sum of the products of the lower right array is $S$. 
\[
  [\lambda_{a}- \lambda_x]_{b}^{l-1} \cdot S 
  + [\lambda_{a}- \lambda_x]_{b+1}^{l-1} \cdot 
  [\lambda_{b}- \lambda_x]_{l}^{L} =
  [\lambda_{a} - \lambda_x]_{b+1}^L 
\] 
 we arrive at
 \[
 [\lambda_{a}- \lambda_x]_{b}^{l-1} S = [\lambda_{a}- \lambda_x]_{b+1}^{l-1}
 \left(
 [\lambda_{a}- \lambda_x]_{l}^{L} -  [\lambda_{b} - \lambda_x]_l^L
 \right)
 \] 
 or equivalently
 \[
 (\lambda_{a}- \lambda_{b}) \cdot S = [\lambda_{a}- \lambda_x]_{l}^{L} -  [\lambda_{b} - \lambda_x]_l^L
 \]
\end{proof}

\begin{lemma}\label{zeroGamma}
For all $1\leq \mu \leq d-2$ we have $\gamma_{\mu,d}=0$. 
\end{lemma}
\begin{proof}
Let $\mu_1=\mu > \mu_2 > \cdots > \mu_s=1 \in \Sigma_\mu$ be a selection of iterations and $d=i_1> i_2 > \cdots \cdots i_s \geq 1> i_{s+1}=0$ be the sequence of $i$'s.  
Using eq. (\ref{eq:lambda-act}) we see that the quantity $[\lambda_{\mu_\nu}^\alpha - \lambda_x]_{i_{\nu+1}+1}^{i_\nu-1}\neq 0$ if and only if one of the following two inequalities hold:
\begin{align}
\label{ineq1}
\text{ either } &&
i_{\nu+1}  >& \mu_\nu - m f(\mu_\nu)
\\
\label{ineq2} 
\text{ or } &&
i_\nu  <&  \mu_\nu+2 -m f(\mu_\nu),
\end{align}
where 
\[
  f(x)=
  \begin{cases}
  1 & \text{ if } m \mid x \\
  0 & \text{ if } m \nmid x
  \end{cases}
\]
We will denote the above two inequalities by (\ref{ineq1})$_\nu$,(\ref{ineq2})$_\nu$ when applied for the integer $\nu$. 
Assume, that for all $1\leq \nu \leq s$ one of the two inequalities (\ref{ineq1})$_\nu$,(\ref{ineq2})$_\nu$ hold, that is $[\lambda_{\mu_\nu}^\alpha - \lambda_x]_{i_{\nu+1}+1}^{i_\nu-1}\neq 0$. 
Inequality (\ref{ineq1})$_{s}$ can not hold for $\nu=s$   
since it gives us $0=i_{s+1}> 1=\mu_s$, we have $m\nmid 1=\mu_s$.

We will keep the sequence $\bar{\mu}:\mu_1> \mu_2 > \cdots > \mu_s$ fixed and we will sum over all possible selections of sequences of $i_1> \cdots i_s > i_{s+1}=0$, that is we will show that the sum 
\begin{equation}
\label{eq:sum-over_i}
\Gamma_{\bar{\mu},i}:=
  \sum_{i=i_1> i_2 > \cdots > i_{s} \geq 1} 
\prod_{\nu=1}^{s}
  [\lambda_{\mu_\nu}^\alpha- \lambda_x]^{i_\nu-1}_{i_{\nu+1}+1}
\end{equation}
is zero, which will show that $\gamma_{\mu,d}=0$ using eq. (\ref{eq:fin-sum4}).

Observe now that if (\ref{ineq2})$_{\nu}$ holds and $m\nmid \nu,\nu-1$, then (\ref{ineq2})$_{\nu-1}$ also holds. Indeed the combination of  (\ref{ineq2})$_{\nu}$  and (\ref{ineq1})$_{\nu-1}$ gives the impossible inequality
\[
\mu_\nu+ 2  \stackrel{(\ref{ineq2})_\nu}{>} i_\nu \stackrel{(\ref{ineq1})_{\nu-1}}{>} \mu_{\nu-1}.
\]
Assume now that $m\mid \nu$ and (\ref{ineq2})$_{\nu}$ holds, then (\ref{ineq2})$_{\nu-1}$ also holds.
Indeed the combination of (\ref{ineq2})$_{\nu}$ and (\ref{ineq1})$_{\nu-1}$
gives us
\[
  \mu_{\nu}+2 -m \stackrel{(\ref{ineq2})_{\nu}}{>} i_{\nu} \stackrel{(\ref{ineq1})_{\nu-1}}{>} \mu_{\nu-1}  - mf(\mu_{\nu-1}). 
\]
If $m \nmid \mu_{\nu-1}$, then the above inequality is impossible since it implies that 
\[
\mu_{\nu} + 2 -m > \mu_{\nu-1} > \mu_{\nu}.
\]
If $m \mid \mu_{\nu-1}$, then the inequality is also impossible since it implies that  $\mu_{\nu}+2 > \mu_{\nu-1}$ so if we write $\mu_{\nu-1}=k' m$ and $\mu_{\nu}= k m$, $k,k' \in \mathbb{N}$, $k'>k$, we arrive at $2>(k'-k)m \geq m$. This proves the following 
\begin{lemma}
\label{when-22}
The inequality (\ref{ineq1})$_{\nu-1}$ might be correct only in cases where 
$m\mid \mu_{\nu-1}$,  $m \nmid \mu_{\nu}$.
\end{lemma}


Assume that for all $\nu$ inequality (\ref{ineq2}) holds. Then for $\nu=1$ it gives us (recall that $\mu \leq d-2$)
\begin{equation}\label{crusial}
  \mu+2 \leq d =i_1 < \mu_1 +2  - mf(\mu_1) = \mu +2  -mf(\mu),
\end{equation}
which is impossible. Therefore either there are $\nu$ such that none of the two inequalities (\ref{ineq1})$_\nu$, (\ref{ineq2})$_\nu$ hold (in this case the contribution to the sum is zero) or there are cases where  (\ref{ineq1}) holds.



The sumands appearing in eq. (\ref{eq:sum-over_i}) can be zero, for example the sequence $\mu_1=m > \mu_2=1$ with $i_2=2< i_1=d$, $s=2$ give the contribution
\[
  [\lambda_{\mu_2}^\alpha - \lambda_x]_1^{i_2-1} [\lambda_{\mu_1}^\alpha-\lambda_x]_{i_2}^{d-1}=
  [\lambda_1^\alpha - \lambda_x]_1^{1} [\lambda_m^\alpha -\lambda_x]^{d-1}_{i_2+1}= (\lambda_2-\lambda_1)
  [\lambda_1 -\lambda_x]^{d-1}_{3}
\]
while for  $i_2=1 < i_1=d$ it gives the contribution 
\[
  [\lambda_{\mu_2}^\alpha - \lambda_x]_1^{i_2-1} [\lambda_{\mu_1}^\alpha-\lambda_x]_{i_2+1}^{d-1}=
  [\lambda_1^\alpha - \lambda_x]_1^{0} [\lambda_m^\alpha -\lambda_x]^{d-1}_{2}= 
  [\lambda_1 -\lambda_x]^{d-1}_{2}
\]
It is clear that these non-zero contributions cancel out when added.


\begin{lemma}
\label{lemma-elim-mi}
Assume that $m\mid \mu_{\nu_0-1}$ and $m \nmid \mu_{\nu_0}$, where 
(\ref{ineq2})$_{\nu_0}$ and (\ref{ineq1})$_{\nu_0-1}$ hold. 
Then, we can eliminate  $\mu_{\nu_0-1}$ and $i_{\nu_0}$ from both selections of the sequence of $\mu$'s and $i$'s, i.e. we can form the sequence of length $s-1$ 
\[
  \bar{\mu}_{s-1}= \mu_s < \bar{\mu}_{s-2}=\mu_{s-1}< \cdots <\bar{\mu}_{\nu_0-1}=\mu_{\nu_0} <
  \bar{\mu}_{\nu_0-2}= \mu_{\nu_0-2}< \cdots  < \bar{\mu}_1 = \mu_1.
\] 
and the corresponding sequence of equal length
\[
\bar{i}_{s-1}= i_s < \bar{i}_{s-2}=i_{s-1}< \cdots <\bar{i}_{\nu_0-1}=i_{\nu_0-1} <
  \bar{i}_{\nu_0}= i_{\nu_0+1}< \cdots  < \bar{i}_1 = i_1=d,
  \]
so that 
\[
\Gamma_{\bar{\mu},i}=
 \sum_{i_1> \cdots > i_s} \prod_{\nu=1}^s [\lambda_{\mu_\nu}^\alpha - \lambda_x]_{i_{\nu+1}+1}^{i_{\nu}-1}
=
(\star)
\sum_{\bar{i}_1> \cdots > \bar{i}_{s-1}} \prod_{\nu=1 \atop \nu\neq \nu_0-1}^s [\lambda_{\mu_\nu}^\alpha - \lambda_x]_{i_{\nu+1}+1}^{i_{\nu}-1},
\]
where $(\star)$ is a non zero element.
\end{lemma}
\begin{proof} (of lemma \ref{lemma-elim-mi})
We are in the case $m\mid \mu_{\nu_0-1}$ and $m \nmid \mu_{\nu_0}$, where 
(\ref{ineq2})$_{\nu_0}$ and (\ref{ineq1})$_{\nu_0-1}$ hold, 
\begin{equation}
\label{ineq3}
  \mu_{\nu_0-1}-m  
    \stackrel{\text{(\ref{ineq1})}_{\nu_0-1}}{<}
  i_{\nu_0} 
  \stackrel{\text{(\ref{ineq2})}_{\nu_0}}{<}
\mu_{\nu_0}+2,
\end{equation} 
or equivalently
\[
\mu_0 := \mu_{\nu_0-1}-m+1 \leq i_{\nu_0} \leq \mu_{\nu_0}+1
\]
For $i_{\nu_0+1}$ the inequality (\ref{ineq1})$_{\nu_0}$  $i_{\nu_0+1} > \mu_{\nu_0} - m f(\mu_{\nu_0})$ can not hold, since it implies 
\[
  i_{\nu_0+1} < i_{\nu_0} \stackrel{(\ref{ineq2})_{\nu_0}}{<} \mu_{\nu_0}+2 < i_{\nu_0+1}+2.
\]


Observe that also 
\[
  i_{\nu_0+1} +1 \leq i_{\nu_0} \leq i_{\nu_{0}-1} -1. 
\]
Set $l=\max\{\mu_0,i_{\nu_0+1}+1\}$ and $L= \min\{\mu_{\nu_0}+1, i_{\nu_0-1}-1\}$. Then $y=i_{\nu_0}$ satisfies 
\[
  l \leq y  \leq L.
\]
By lemma \ref{lemma:ABlL} the quantity
\[
 \sum_{l \leq y \leq L} 
 [\lambda_{\mu_{\nu_0}+1} - \lambda_x]_{i_{\nu_0+1}+1}^{y-1}
 \cdot
 [\lambda_{\mu_0} - \lambda_x]_{y+1}^{i_{\nu_0-1}-1} 
\]
equals to
\[
[\lambda_{\mu_{\nu_0}+1} - \lambda_x]_{i_{\nu_0+1}+1}^{l-1}
 \cdot
 [\lambda_{\mu_0} - \lambda_x]_{L+1}^{i_{\nu_0-1}-1}
  \cdot 
\frac{
  [\lambda_{\mu_{\nu_0}+1}- \lambda_x]_{l}^{L} -  [\lambda_{\mu_0} - \lambda_x]_l^L
}
{(\lambda_{\mu_{\nu_0}+1}- \lambda_{\mu_0})}
\]
\begin{equation}
\label{eq:after-sum}
  \frac{
  [\lambda_{\mu_{\nu_0}+1}- \lambda_x]_{i_{\nu_0+1}+1}^{L}
  \cdot
  [\lambda_{\mu_0} - \lambda_x]_{L+1}^{i_{\nu_0-1}-1}
   -  
  [\lambda_{\mu_{\nu_0}+1} - \lambda_x]_{i_{\nu_0+1}+1}^{l-1}\cdot 
  [\lambda_{\mu_0} - \lambda_x]_l^{i_{\nu_0-1}-1}
}
{(\lambda_{\mu_{\nu_0}+1}- \lambda_{\mu_0})}
\end{equation}

\noindent {\bf Case A1} $l=\mu_0 \geq i_{\nu_0+1}+1$. Then $[\lambda_{\mu_0} - \lambda_x]_l^L=0$.

\noindent {\bf Case A2} $l=i_{\nu_0+1}+1 > \mu_0$. We set $z:=i_{\nu_0+1}$, which is bounded by eq. (\ref{ineq2})$_{\nu_{0}+1}$ that is 
\[
   \mu_0 \stackrel{\text{Case A2}}{\leq}  z  \stackrel{(\ref{ineq2})_{\nu_0+1}}{\leq} \mu_{\nu_0+1}+1.
\]
Notice that in this case $m\nmid \mu_{\nu_0+1}$. Indeed, we have assumed that inequality (\ref{ineq2})$_{\nu_0+1}$ holds wich gives us 
\[
\mu_{\nu_0-1}-m =\mu_0-1 \stackrel{(\text{Case A2})}{<} i_{\nu_0+1} 
\stackrel{(\ref{ineq2})_{\nu_0+1}}{<} \mu_{\nu_0+1}+2 -m,
\]
which implies that $\mu_{\nu_0-1} < \mu_{\nu_0+1}+2$, a contradiction. 
Thus for $l=z+1$ we compute
\[
  \sum_{ \mu_0 \leq z \leq \mu_{\nu_0+1}+1}
  [\lambda_{\mu_{\nu_0+1}}^\alpha - \lambda_x]^{i_{\nu_0+1}-1}_{i_{\nu_0+2}+1} 
  \cdot
  [\lambda_{\mu_0} -\lambda_x]^L_l
  =
\]
\[
  =\sum_{\mu_0 \leq z \leq \mu_{\nu_0+1}+1}
  [\lambda_{\mu_{\nu_0+1}+1} - \lambda_x]^{z-1}_{i_{\nu_0+2}+1} 
  \cdot
  [\lambda_{\mu_0} -\lambda_x ]^L_{z+1}=
\]
\[
  =
  (\star) \cdot
  \frac{
  [\lambda_{\mu_{\nu_0+1}+1}- \lambda_x]_{\mu_0}^{\mu_{\nu_0+1}+1}
  -
  [\lambda_{\mu_0} - \lambda_x]_{\mu_0}^{\mu_{\nu_0+1}+1}
  }
  {\lambda_{\mu_{\nu_0+1}+1}-\lambda_{\mu_{0} +1 }}=0.
\]
\noindent{\bf Case B1} $L=\mu_{\nu_0}+1 \leq  i_{\nu_0-1}-1$. In this case $[\lambda_{\mu_{\nu_0}+1} - \lambda_x]_l^L=0$. 

\noindent{\bf Case B2} $L=i_{\nu_0-1}-1 < \mu_{\nu_0}+1$. 
In this case eq. (\ref{eq:after-sum}) is reduced to 
\[
 \frac{
 [\lambda_{\mu_{\nu_0}+1} - \lambda_x]_{i_{\nu_0+1}+1}^{i_{\nu_0-1}-1}
 }
 {
 (\lambda_{\mu_{\nu_0}+1}- \lambda_{\mu_0})
 }
\]
This means that we have erased the $\mu_{\nu_0-1}$ from the product and we have
\[
 \sum_{i_1> \cdots > i_s} \prod_{\nu=1}^s [\lambda_{\mu_\nu}^\alpha - \lambda_x]_{i_{\nu+1}+1}^{i_{\nu}-1}
=
(\star)
\sum_{i_1> \cdots > i_s} \prod_{\nu=1 \atop \nu\neq \nu_0-1}^s [\lambda_{\mu_\nu}^\alpha - \lambda_x]_{i_{\nu+1}+1}^{i_{\nu}-1},
\]
where $(\star)$ is a non zero element.
This procedure gives us that the original quantity 
\[
[\lambda_{\mu_{\nu_0}}^\alpha - \lambda_x]_{i_{\nu_0+1}+1}^{i_{\nu_0}-1} \cdot 
[\lambda_{\mu_{\nu_0-1}}^\alpha -\lambda_x]_{i_{\nu_0}+1}^{i_{\nu_0-1}-1}
\]
after summing over $i_{\nu_0}$ becomes the quantity
\[
[\lambda_{\mu_{\nu_0}}^\alpha - \lambda_x]_{i_{\nu_0+1}+1}^{i_{\nu_0-1}-1} = 
[\lambda_{\bar{\mu}_{\nu_0-1}}^\alpha - \lambda_x]_{\bar{i}_{\nu_0}+1}^{\bar{i}_{\nu_0-1}-1},
\]
that is we have eliminated the $\mu_{\nu_0-1}$ and $i_{\nu_0}$ from both selections of the sequence of $\mu$'s and $i$'s, i.e. we have the sequence of length $s-1$ 
\[
  \bar{\mu}_{s-1}= \mu_s < \bar{\mu}_{s-2}=\mu_{s-1}< \cdots <\bar{\mu}_{\nu_0-1}=\mu_{\nu_0} <
  \bar{\mu}_{\nu_0-2}= \mu_{\nu_0-2}< \cdots  < \bar{\mu}_1 = \mu_1.
\] 
and the corresponding sequence of equal length
\[
\bar{i}_{s-1}= i_s < \bar{i}_{s-2}=i_{s-1}< \cdots <\bar{i}_{\nu_0-1}=i_{\nu_0-1} <
  \bar{i}_{\nu_0}= i_{\nu_0+1}< \cdots  < \bar{i}_1 = i_1=d.
\] 
\end{proof}

\begin{remark}
\label{strong-inequality}
One should be careful here since  $ \bar{i}_{\nu_0-1}=i_{\nu_0-1}> i_{\nu_0}  >\bar{i}_{\nu_0}=i_{\nu_0+1} $, so $\bar{i}_{\nu_0-1}> \bar{i}_{\nu_0}+1$. This means that the new sequence of $\bar{i}_{s-1}> \cdots > \bar{i}_1$ satisfies a stronger inequality in the $\nu_0$ position, unless $\nu_0-1=d$ in the computation of $\gamma_{d,d}$.  
\end{remark}

Consider the set $s,s-1,\ldots,\nu_0$ such that $m\nmid \mu_{\nu}$ for $s \geq \nu \geq \nu_0$ and assume that $m \mid \mu_{\nu_0-1}$ and (\ref{ineq2})$_{\nu_0}$ and (\ref{ineq1})$_{\nu_0-1}$ hold. 
 We apply lemma 
\ref{lemma-elim-mi} and we obtain a new sequence of $\mu$'s with $\mu_{\nu_0-1}$ removed, provided that $\nu_0-1>1$. 
We continue this way and in the sequence of $\mu$'s we eliminate all possible inequalities like (\ref{ineq3})
obtaining a series of $\mu$ which involves only inequalities of type (\ref{ineq2}). But this is not possible if $\mu \leq  d-2$, according to equation (\ref{crusial}). This proves that all $\gamma_{\mu,d}=0$ for $1\leq \mu \leq d-2$, this completes the proof of lemma \ref{zeroGamma}. 
\end{proof}

\begin{lemma}\label{d-1}
If $\mu_2\neq d-1$, then the  contribution of the corresponding summand $\Gamma_{\bar{\mu},i}$ to  $\gamma_{d,d}$ is zero. 
\end{lemma}
\begin{proof}
We are in the case $\mu=d=i$. 
We begin the procedure of eliminating all sequences of inequalities of the form $(23)_{\nu_0},(22)_{\nu_0-1}$, where $m\mid \nu_0-1$, $m \nmid \nu_0$, using lemma \ref{lemma-elim-mi}. 
For $\nu=1$ inequality (\ref{ineq2})$_1$ can not hold since it implies the impossible inequality $d=i_1 < d+2 -m$. Therefore, (\ref{ineq1})$_1$ holds, that is $i_2> d -m$. 
On the other hand we can assume that (\ref{ineq2})$_2$ holds by the elimination process,  so we have
\[
  d-m \stackrel{(\ref{ineq1})_1}{<} i_2 \stackrel{(\ref{ineq2})_2} < \mu_2+2.
\]
Following the analysis of the proof of lemma \ref{zeroGamma} we see that 
the contribution to 
$\gamma_{d,d}$ is non zero if case B2 holds, that is ($\nu_0=2$ in this case)  $d-1=i_{\nu_0-1}-1 < \mu_2+1$, obtaining that  $\mu_2=d-1$. 
\end{proof}
\begin{lemma}\label{lemma:30}
Equation (\ref{finalEquation}) holds, that is 
\[
   (\lambda_d - \lambda_d^{\alpha}) \gamma_{d,d} =  \sum_{\nu=1}^{d-1} t_{d,\nu}^{(\alpha)} \gamma_{\nu,d}=t_{d,d-1}^{(\alpha)} \gamma_{d-1,d}.
\]
\end{lemma}
\begin{proof}
We will use the procedure of the proof of lemma \ref{lemma-elim-mi}. We recall that for each fixed sequence of $\mu_s >\cdots > \mu_1$ we summed over all possible sequences $i_1 > \cdots  >i_{s+1}=0$. 
In the final step the inequality (\ref{ineq3}) appears, for $\nu_0=2$, and $\mu_{\nu_0}=\mu_{2}=d-1$ and $\nu_0-1=1$ and $\mu_{\nu_0-1}=\mu=d$, that is: 
\[
   0=\mu_{\nu_0-1}-m  
    \stackrel{\text{(\ref{ineq1})}_{2}}{<}
  i_{\nu_0} 
  \stackrel{\text{(\ref{ineq2})}_{1}}{<}
\mu_{\nu_0}+2=d+1.
\]
As in the proof of lemma \ref{lemma-elim-mi} we sum over $y=i_{\nu}$ and the result is either zero in case B1 or in the B2 case, where  $\mu_{\nu_0}=\mu_2=d-1$ and $\mu_{0}=\mu_{\nu_0-1}-m+1=d-m+1$, the contribution is computed to be equal to
\[
  \frac{
[\lambda_{\mu_{\nu_0}+1}^\alpha - \lambda_x]_{i_{\nu_0+1}+1}^{i_{\nu_0-1}-1}
  }{
(\lambda_{\mu_{\nu_0}+1} - \lambda_{m_0})
  }
  =
  \frac{
  [\lambda_{d}^\alpha - \lambda_x]_{i_{3}+1}^{d-1}
  }{\lambda_d - \lambda_d^\alpha}.
\]
The last $\mu_{\nu_0-1}=\mu_1=d$ is eliminated in the above expression. 
This means that for a fixed sequence $\mu_1 > \ldots > \mu_s$ the contribution of the inner sum in eq. (\ref{eq:sum-over_i}) is given by
\[
\frac{1}{\lambda_d- \lambda_d^\alpha} \cdot
  \sum_{d-1=i_2> i_3 > \cdots > i_{s} \geq 1} 
\prod_{\nu=2}^{s}
  [\lambda_{\mu_\nu}^\alpha- \lambda_x]^{i_\nu-1}_{i_{\nu+1}+1}.
\]
Observe that $\mu_1=d$ does not appear in this expression and this expression corresponds to the sequence $\bar{\mu}_1=\mu_2=d-1 > \bar{\mu}_2=\mu_3 > \cdots > \bar{\mu}_{s-1}=\bar{\mu}_s=1$.
Notice, also that the problem described in remark \ref{strong-inequality} does not appear here, sence we erased  $i_1$ which is not between some $i$'s but the first one.   
Therefore, 
we can relate it to a similar expression that contributes to $\gamma_{d-1,d}$. Conversely every contribution of $\gamma_{d-1,d}$ gives rise to a contribution in $\gamma_{d,d}$, by multiplying by 
$\lambda_d- \lambda_d^\alpha$. The desired result follows by the expression of $\gamma_{\mu,d}$ given in eq. (\ref{eq:fin-sum4}).  
\end{proof}

We have shown so far how to construct matrices $\Gamma, T$ so that 
\begin{equation}
T^q=1, \Gamma T \Gamma^{-1}=T^\alpha.
\label{eq:partial-rel}
\end{equation} 
We will now prove that $\Gamma$ has order $m$. By equation (\ref{eq:partial-rel}) $\Gamma^k$ should satisfy equation 
\[
  \Gamma^k T \Gamma^{-k} =T^{\alpha^k}.
\]
Using proposition \ref{prop:fing-G-uniqe} asserting the uniqueness of such $\Gamma^k$ with $\alpha$ replaced by $\alpha^k$ we have that the matrix multiplication of the entries of $\Gamma$ giving rise to $(\gamma_{\mu,i}^{(k)})=\Gamma^k$ coincide to the values by the the recursive method of proposition (\ref{eq:partial-rel}) applied for $\Gamma'=\Gamma^k$, $\alpha'=\alpha^k$ and $\Gamma'E_1= \zeta_m^{\epsilon k} E_1$. In particular for $k=m$, we have $\alpha^m\equiv 1 \mod p^\nu$ for all $1\leq \nu \leq h$, that is the matrix $\Gamma^k$ should be recursively constructed using proposition (\ref{eq:partial-rel}) for the relation $\Gamma^m T \Gamma^m =T$, $\Gamma^m E_1=E_1$, leading to the conclusion $\Gamma^m=\mathrm{Id}$. Notice that the first eigenvalue of $\Gamma$ is a primitive root of unity, therefore $\Gamma$ has order exactly $m$. 

By lemma \ref{action-tau} the action of $\sigma$ in the special fibre is given by a lower triangular matrix. Therefore, we must have
\begin{equation}
\label{eq:condition-lower-triangular}
 \gamma_{\nu,i} \in \mathfrak{m}_r \text{ for } \nu < i.
\end{equation}
\begin{proposition}
\label{prop:25}
If 
\begin{equation}
\label{eq:divis-cond}
v(\lambda_i-\lambda_j) > v(a_\nu) \text{ for all } 1\leq i,j \leq d  \text{ and } 1\leq \nu \leq d-1
\end{equation}
 then the 
matrix $(\gamma_{\mu,i})$ has entries in the ring $R$ and is lower triangular modulo $\mathfrak{m}_R$.
\end{proposition}
\begin{proof}
Assume that the condition of eq. (\ref{eq:divis-cond}) holds. 
In equation (\ref{eq:fin-sum4}) we compute the fraction
\begin{equation}
\label{eq:aaa}
\frac{[a]_1^{\mu-1}}{[a]_1^{i-1}}=
\begin{cases}
\frac{1}{[a]_{\mu}^{i-1}
} & \text{ if } i > \mu 
\\ 
1 & \text{ if } i=\mu
\\
[a]^{\mu-1}_{i} & \text{ if } i < \mu 
\end{cases}
\end{equation}
The number of $(\lambda_\mu^\alpha- \lambda_x)$ factors in the numerator is equal to (recall that $i_{s+1}=0$)
\[
  \sum_{\nu=1}^{s} (i_\nu-1- i_{\nu+1}-1+1) 
  = i -s,
\]
and $i > \mu \geq s$, so  $i-s>0$. 
Therefore, for the upper part of the matrix $i > \mu$ we have $i-s$  factors of the form $(\lambda_i^\alpha-\lambda_j)$
in the numerator and $i-\mu$ factors $a_x$ in the denominator. Their difference is equal to $(i-s)-(i-\mu)= \mu- s \geq 0$. 
By assumption the matrix reduces to an upper triangular matrix modulo $\mathfrak{m}_R$.
\end{proof}

\begin{remark}
\label{cond-ram-R}
The condition given in equation (\ref{eq:divis-cond}) can be satisfied in the following way: It is clear that $\lambda_i-\lambda_j \in \mathfrak{m}_R$. Even in the case $v_{\mathfrak{m}_R}(\lambda_i - \lambda_j)=1$ we can consider a ramified extension $R'$ of the ring $R$ with ramification index $e$, in order to make the valuation $v_{\mathfrak{m}_{R'}}(\lambda_i-\lambda_j)=e$ and then there is space to select $v_{\mathfrak{m}_{R'}}(a_i)< v_{\mathfrak{m}_{R'}}(\lambda_i-\lambda_j)$. 
\end{remark}


\begin{proposition}
\label{prop:reduction}
We have that 
\begin{equation}
\label{eq:E1-detrmines-all}
  \gamma_{i,i}\equiv \zeta_m^\epsilon  \alpha^{i-1} \mod \mathfrak{m}_R
\end{equation}
Let $A=\{a_1,\ldots,a_{d-1}\} \in R$ be the set of elements below the diagonal in eq. (\ref{eq-matrix-lift}).  If $a_i\in \mathfrak{m}_R$, then 
\[
\gamma_{\mu,i} \in \mathfrak{m}_R \text{ for } \mu\neq i,
\]
that is $E_i$ is an eigenvector for the reduced action of $\Gamma$ modulo $\mathfrak{m}_R$. 
If $a_{\kappa_1},\ldots,a_{\kappa_r}$ the elements of the set $A$ which are in $\mathfrak{m}_R$, then the reduced matrix of $\Gamma$ has the form:
\[
\begin{pmatrix}
\Gamma_1 & 0 & \cdots &  0 
\\
0 & \Gamma_2        & \ddots & \vdots 
\\  
\vdots & \ddots & \ddots & 0
\\ 
0 & \cdots & 0 & \Gamma_r 
\end{pmatrix}
\]
where $\Gamma_1,\Gamma_2,\ldots,\Gamma_{r+1}$ for $1\leq \nu \leq r+1$ are $(\kappa_\nu-\kappa_{\nu-1}) \times (\kappa_\nu-\kappa_{\nu-1})$  lower triangular matrices (we set $\kappa_0=0, \kappa_{r+1}=d$).
\end{proposition}
\begin{proof}
Consider the matrix $\Gamma$:
\[
\begin{pmatrix}
\tikzmarknode{1A}{
\gamma_{11} 
}&  & &  
\\
\vdots & \ddots & & &   & & & \bigzero  \\
\;\;\gamma_{\kappa_1,1} &  \cdots &
\tikzmarknode{1B}{
  \gamma_{\kappa_1,\kappa_1}
  }
\\
\tikzmarknode{1C}{\color{white}\gamma_{11}}
& & & 
\tikzmarknode{2A}{
\gamma_{\kappa_1+1,\kappa_1+1}
} & & 
\\
& & & \vdots & \ddots &  & &   \\
& 
\gamma_{\mu,i} &
 & \gamma_{\kappa_2,\kappa_1+1} &  \cdots &
\tikzmarknode{2B}{
\gamma_{\kappa_2,\kappa_2}
}
\\ 
 & \tikzmarknode{MM} & &
\tikzmarknode{2C}{\color{white}\gamma_{\kappa_1+1,\kappa_1+1}}
 & \gamma_{\mu,i} & & \ddots 
\\
&   && &   \tikzmarknode{MM1} & 
& 
&  
\tikzmarknode{3A}{
\gamma_{\kappa_r+1,\kappa_r+1}
} 
\\
&&  &  &  & & \cdots & \vdots & \ddots &   \\
&& \tikzmarknode{1D}{
 {\color{white}{\gamma_{\kappa_1,\kappa_1}}} 
}& & &   
\tikzmarknode{2D}{
 {\color{white}{\gamma_{\kappa_2,\kappa_2}}}}
 & 
& \gamma_{d,\kappa_r+1} &  \cdots & 
\tikzmarknode{3B}{
  \!\!\gamma_{d,d}  
  } 
\end{pmatrix}
\]
\begin{tikzpicture}[overlay, remember picture]
\foreach \X in {1,2,3}
{\node[inner sep=0.9pt,draw,color=red,fit=(\X A)(\X B)]{};}
\foreach \X in {1,2}
{\node[inner sep=0.9pt,draw,color=blue,fit=(\X C)(\X D)]{};}
\node at (MM) { \small $\;\;1\leq i \leq \kappa_1 < m \leq d$};
\node at (MM1) { \small $\!\!\!\!\!\!\!\!\kappa_1 < i \leq \kappa_2<\mu \leq d$};
\end{tikzpicture}
We have that $\mu=i$ and the only element in $\Sigma_\mu$ which does not have any factor of the form $(\lambda_y^\alpha - \lambda_x )$ is the sequence 
\[
1=\mu_s=\mu_{s-1}-1 < \mu_{s-1}< \cdots <   \mu_2=\mu_1-1 <\mu_1=\mu
\]
For this sequence eq. (\ref{eq:fin-sum4}) becomes
\begin{equation*}
 \gamma_{i,i}= \prod_{\nu=2}^{s} h_{\alpha-1}(\lambda_{\mu_\nu},\lambda_{\mu_{\nu-1}}) \zeta_{m}^\epsilon \mod \mathfrak{m}_R,
\end{equation*}
which gives the desired result since $h_{\alpha-1}(\lambda_{\mu_\nu},\lambda_{\mu_{\nu-1}})\equiv \binom{\alpha}{1}= \alpha \mod \mathfrak{m}_R$.

For proving that all entries $\gamma_{\mu,i} \in \mathfrak{m}_R$ 
for $\kappa_\nu < i \leq \kappa_{\nu+1 } < \mu \leq d $, that is for all entries bellow the central blocks,  we observe that 
from equation (\ref{eq:fin-sum4}) combined with eq. (\ref{eq:aaa}) that $\gamma_{\mu,i}$ is divisible by $[a]_i^{\mu-1}=a_i a_{i+1} \cdots a_{\kappa_\nu+1} \cdots  a_{\mu-1} \in \mathfrak{m}_R$.



\end{proof}

Recall that by lemma \ref{loga} there is an $1\leq a_0 \leq m$ such that $\alpha=\zeta_m^{a_0}$. 
\begin{proposition}
The indecomposable module $V$ modulo $\mathfrak{m}_R$ breaks into a direct sum of $r+1$ indecomposable $k[C_q \rtimes C_m]$ modules $V_{\nu}$, $1\leq \nu \leq r+1$. Each $V_\nu$ is isomorphic to $V_{\alpha}(\epsilon+a_0 \kappa_{\nu-1},\kappa_\nu-\kappa_{\nu-1})$.
\end{proposition}
\begin{proof}
By eq. (\ref{eq:E1-detrmines-all}) the first eigenvalue of the reduced matrix block $\Gamma_\nu$
is 
\[
\zeta_m^{\epsilon} \alpha^{\kappa_{\nu-1}}=\zeta_m^{\epsilon+ (\kappa_{\nu-1})a_0}.
\] 
Since that first eigenvalue together with the size of the block determine the last eigenvalue, that is the action of $C_m$ on the socle the reduced block is uniquely determined up to isomorphism. 
\end{proof}

This way we arrive at a new obstruction. Assume that the indecomposable representation given by the matrix $T$ as in lemma 
\ref{lemma:Jordan} reduces modulo $\mathfrak{m}_R$ to a sum of Jordan blocks. Then the $\sigma$ action on the leading elements of each Jordan block in the special fibre should be described by the corresponding action of $\sigma$ on the leading eigenvector $E$ of $V$. The corresponding actions on the special fibre should be compatible. 


This observation is formally given in proposition \ref{cond-lift}, which we now prove: 
Each set $I_\nu$, $1\leq \nu \leq t$ corresponds to an indecomposable $R[G]$-module, which decomposes to the indecomposables $V_{\alpha}(\epsilon_\mu,\kappa_\mu)$, $\nu\in I_\nu$ of the special fiber. Indecomposable summands have different roots of unity in $R$, therefore $\sum_{\mu\in I_\nu} k_\nu \leq q$, this is condition (\ref{cond-lift}.\ref{cond-lifta}). The second condition  (\ref{cond-lift}.\ref{cond-liftb}) comes from proposition \ref{prop:div-obstr}. If $1$ is one of the possible eigenvalues of the lift $T$, then 
$\sum_{\mu \in I_\nu} \kappa_\mu \equiv 1 \mod m$. If all eigenvalues of the lift $T$ are different than one, then $\sum_{\mu \in I_\nu} \kappa_\mu \equiv 0 \mod m$. If $\#I_\nu=q$, then there is one zero eigenvalue and the sum equals $1\mod m$.

It is clear by eq. (\ref{eq:E1-detrmines-all}) that condition (\ref{cond-lift}.\ref{cond-liftc})  is a necessary condition. 
On the other hand if (\ref{cond-lift}.\ref{cond-liftc}) is satisfied we can write (after a permutation if necessary) the set $\{1,\ldots,S\}$, $S=\sum_{\nu=1}^t \#I_\nu$ as
\begin{align*}
J_1 &=\{1,2,\ldots,\kappa_1^{(1)},\kappa_1^{(1)}+1,\ldots,\kappa_1^{(1)}+\kappa_2^{(1)},\ldots,
\sum_{j=1}^{r_1} \kappa_j^{(1)}=b_1\},  I_1=\{\kappa_1^{(1)},\ldots,\kappa_{r_1}^{(1)}\}
\\
J_2 &= \{b_1+1, b_1+2,\ldots,b_2=b_1+\sum_{j=1}^{r_2} \kappa_j^{(2)}\}, I_2=\{\kappa_1^{(2)},\ldots,\kappa_{r_2}^{(2)}\}
\\
\cdots & \cdots 
\\
J_s &= \{b_{s-1}+1,b_{t-1}+2,\ldots, b_t=S\}, I_s=\{\kappa_1^{(s)},\ldots,\kappa_{r_s}^{(s)}\}
\end{align*}
The matrix given in eq. (\ref{eq-matrix-lift}), where 
\[
a_i=
\begin{cases}
0 & \text{ if } i \in \{b_1,\ldots,b_{s-1}\} \\
\pi & \text{ if } i \in 
\{
\kappa_{1}^{(\nu)},
\kappa^{(\nu)}_{1}+\kappa^{(\nu)}_{2},
\kappa^{(\nu)}_{1}+\kappa^{(\nu)}_{2}+\kappa^{(\nu)}_{3},
\ldots,
\kappa^{(\nu)}_{1}+\kappa^{(\nu)}_{2}+\cdots+\kappa^{(\nu)}_{r_{\nu}-1}\} \\
1 & \text{ otherwise}  
\end{cases}  
\] 
lifts the $\tau$ generator, and by (\ref{eq:allcombined}) there is a well defined extended action of the $\sigma$ as well. 

\noindent {\bf Example:}
Consider the group $q=5^2,m=4,\alpha=7$,
\[
G=C_{5^2} \rtimes C_4= 
\langle
\sigma,\tau | \sigma^4= \tau^{25}=1, \sigma \tau \sigma^{-1} = \tau^{7} 
\rangle.
\]
Observe that $\mathrm{ord}_5 7= \mathrm{ord}_{5^2} 7=4$. 
\begin{itemize}
  \item
The module $V(\epsilon,25)$ is projective and is known to lift in characteristic zero. This fits well with proposition \ref{cond-lift}, since $4 \mid 25-1=4\cdot 6$. 
 \item 
The modules $V(\epsilon, \kappa)$ do not lift in characteristic zero if $4\nmid \kappa$ or $4 \nmid \kappa-1$. Therefore only $V(\epsilon,1)$, $V(\epsilon,4)$, $V(\epsilon,5)$, $V(\epsilon,8)$, $V(\epsilon,9)$, $V(\epsilon,12)$, $V(\epsilon,13)$, $V(\epsilon,16)$, $V(\epsilon,17)$, $V(\epsilon,20)$, $V(\epsilon,21)$, $V(\epsilon,24)$, $V(\epsilon,25)$ lift.
\item 
The module $V(1,2)\oplus V(3,2)$ lift to characteristic zero, where the matrix of $T$ with respect to a basis $E_1,E_2,E_3,E_4$ is given by  
\[
 T=
\begin{pmatrix}
\zeta_q & 0 & 0 & 0 \\
1 & \zeta_q^2 & 0 & 0 \\
0 & \pi &  \zeta_q^3 & 0 \\
0 &  0 &    1 & \zeta_q^4  
\end{pmatrix} 
\]
and $\sigma(E_1)= \zeta_q E_1$. 
\item 
The module $V(1,2) \oplus V(1,2)$ does not lift in characteristic zero. There is no way to permute the direct summands so that the eigenvalues of $\sigma$ are given by $\zeta_m^{\epsilon}, \alpha \zeta_m^{\epsilon}, \alpha^2\zeta_m^{\epsilon}, \alpha^3 \zeta_m^{\epsilon}$. Notice that $\alpha=2=\zeta_m$. 
\item 
The module $V(\epsilon_1,21)\oplus V(2^{21}\cdot \epsilon_1,23)$ does not lift in characteristic zero. The sum $21+24$ is divisible by $4$, $\epsilon_2=2^{21} \epsilon_1$ is compatible, but $21+23=44>25$ so the representation of $T$ in the supposed indecomposable module formed by their sum can not have different eigenvalues which should be $25$-th roots of unity. 
\end{itemize}

 \def\cprime{$'$}

\end{document}